\documentstyle[amssymb,amstex,amscd,twoside,leqno]{amsart}
\makeatletter
 
\def\diagram{\m@th\leftwidth=\z@ \rightwidth=\z@ \topheight=\z@
\botheight=\z@ \setbox\@picbox\hbox\bgroup}
 
\def\enddiagram{\egroup\wd\@picbox\rightwidth\unitlength
\ht\@picbox\topheight\unitlength \dp\@picbox\botheight\unitlength
\hskip\leftwidth\unitlength\box\@picbox}
 
\def\bfig{\begin{diagram}}
\def\efig{\end{diagram}}
\newcount\wideness \newcount\leftwidth \newcount\rightwidth
\newcount\highness \newcount\topheight \newcount\botheight
 
\def\ratchet#1#2{\ifnum#1<#2 \global #1=#2 \fi}
 
\def\putbox(#1,#2)#3{%
\horsize{\wideness}{#3} \divide\wideness by 2
{\advance\wideness by #1 \ratchet{\rightwidth}{\wideness}}
{\advance\wideness by -#1 \ratchet{\leftwidth}{\wideness}}
\vertsize{\highness}{#3} \divide\highness by 2
{\advance\highness by #2 \ratchet{\topheight}{\highness}}
{\advance\highness by -#2 \ratchet{\botheight}{\highness}}
\put(#1,#2){\makebox(0,0){$#3$}}}
 
\def\putlbox(#1,#2)#3{%
\horsize{\wideness}{#3}
{\advance\wideness by #1 \ratchet{\rightwidth}{\wideness}}
{\ratchet{\leftwidth}{-#1}}
\vertsize{\highness}{#3} \divide\highness by 2
{\advance\highness by #2 \ratchet{\topheight}{\highness}}
{\advance\highness by -#2 \ratchet{\botheight}{\highness}}
\put(#1,#2){\makebox(0,0)[l]{$#3$}}}
 
\def\putrbox(#1,#2)#3{%
\horsize{\wideness}{#3}
{\ratchet{\rightwidth}{#1}}
{\advance\wideness by -#1 \ratchet{\leftwidth}{\wideness}}
\vertsize{\highness}{#3} \divide\highness by 2
{\advance\highness by #2 \ratchet{\topheight}{\highness}}
{\advance\highness by -#2 \ratchet{\botheight}{\highness}}
\put(#1,#2){\makebox(0,0)[r]{$#3$}}}

\def\adjust[#1]{} 
 
\newcount \coefa
\newcount \coefb
\newcount \coefc
\newcount\tempcounta
\newcount\tempcountb
\newcount\tempcountc
\newcount\tempcountd
\newcount\xext
\newcount\yext
\newcount\xoff
\newcount\yoff
\newcount\gap%
\newcount\arrowtypea
\newcount\arrowtypeb
\newcount\arrowtypec
\newcount\arrowtyped
\newcount\arrowtypee
\newcount\height
\newcount\width
\newcount\xpos
\newcount\ypos
\newcount\run
\newcount\rise
\newcount\arrowlength
\newcount\halflength
\newcount\arrowtype
\newdimen\tempdimen
\newdimen\xlen
\newdimen\ylen
\newsavebox{\tempboxa}%
\newsavebox{\tempboxb}%
\newsavebox{\tempboxc}%
 
\newdimen\w@dth
 
\def\setw@dth#1#2{\setbox\z@\hbox{\m@th$#1$}\w@dth=\wd\z@
\setbox\@ne\hbox{\m@th$#2$}\ifnum\w@dth<\wd\@ne \w@dth=\wd\@ne \fi
\advance\w@dth by 1.2em}
 
\def\t@^#1_#2{\allowbreak\def\n@one{#1}\def\n@two{#2}\mathrel
{\setw@dth{#1}{#2}
\mathop{\hbox to \w@dth{\rightarrowfill}}\limits
\ifx\n@one\empty\else ^{\box\z@}\fi
\ifx\n@two\empty\else _{\box\@ne}\fi}}
\def\t@@^#1{\@ifnextchar_{\t@^{#1}}{\t@^{#1}_{}}}
\def\to{\@ifnextchar^{\t@@}{\t@@^{}}}
 
\def\t@left^#1_#2{\def\n@one{#1}\def\n@two{#2}\mathrel{\setw@dth{#1}{#2}
\mathop{\hbox to \w@dth{\leftarrowfill}}\limits
\ifx\n@one\empty\else ^{\box\z@}\fi
\ifx\n@two\empty\else _{\box\@ne}\fi}}
\def\t@@left^#1{\@ifnextchar_{\t@left^{#1}}{\t@left^{#1}_{}}}
\def\toleft{\@ifnextchar^{\t@@left}{\t@@left^{}}}
 
\def\two@^#1_#2{\allowbreak
\def\n@one{#1}\def\n@two{#2}\mathrel{\setw@dth{#1}{#2}
\mathop{\vcenter{\lineskip\z@\baselineskip\z@
                 \hbox to \w@dth{\rightarrowfill}%
                 \hbox to \w@dth{\rightarrowfill}}%
       }\limits
\ifx\n@one\empty\else ^{\box\z@}\fi
\ifx\n@two\empty\else _{\box\@ne}\fi}}
\def\tw@@^#1{\@ifnextchar _{\two@^{#1}}{\two@^{#1}_{}}}
\def\two{\@ifnextchar ^{\tw@@}{\tw@@^{}}}
 
\def\tofr@^#1_#2{\def\n@one{#1}\def\n@two{#2}\mathrel{\setw@dth{#1}{#2}
\mathop{\vcenter{\hbox to \w@dth{\rightarrowfill}\kern-1.7ex
                 \hbox to \w@dth{\leftarrowfill}}%
       }\limits
\ifx\n@one\empty\else ^{\box\z@}\fi
\ifx\n@two\empty\else _{\box\@ne}\fi}}
\def\t@fr@^#1{\@ifnextchar_ {\tofr@^{#1}}{\tofr@^{#1}_{}}}
\def\tofro{\@ifnextchar^ {\t@fr@}{\t@fr@^{}}}

\def\mon{\mathop{\m@th\hbox to
      14.6\P@{\lasyb\char'51\hskip-2.1\P@$\arrext$\hss
$\mathord\rightarrow$}}\limits} 
\def\leftmono{\mathrel{\m@th\hbox to
14.6\P@{$\mathord\leftarrow$\hss$\arrext$\hskip-2.1\P@\lasyb\char'50%
}}\limits} 
\mathchardef\arrext="0200       

\def\settypes(#1,#2,#3){\arrowtypea#1 \arrowtypeb#2 \arrowtypec#3}
\def\settoheight#1#2{\setbox\@tempboxa\hbox{#2}#1\ht\@tempboxa\relax}%
\def\settodepth#1#2{\setbox\@tempboxa\hbox{#2}#1\dp\@tempboxa\relax}%
\def\settokens`#1`#2`#3`#4`{%
     \def\tokena{#1}\def\tokenb{#2}\def\tokenc{#3}\def\tokend{#4}}
\def\setsqparms[#1`#2`#3`#4;#5`#6]{%
\arrowtypea #1
\arrowtypeb #2
\arrowtypec #3
\arrowtyped #4
\width #5
\height #6
}
\def\setpos(#1,#2){\xpos=#1 \ypos#2}

\def\settriparms[#1`#2`#3;#4]{\settripairparms[#1`#2`#3`1`1;#4]}%
 
\def\settripairparms[#1`#2`#3`#4`#5;#6]{%
\arrowtypea #1
\arrowtypeb #2
\arrowtypec #3
\arrowtyped #4
\arrowtypee #5
\width #6
\height #6
}
 
\def\resetparms{\settripairparms[1`1`1`1`1;500]\width 500}
 
\resetparms
 
\def\mvector(#1,#2)#3{
\put(0,0){\vector(#1,#2){#3}}%
\put(0,0){\vector(#1,#2){26}}%
}
\def\evector(#1,#2)#3{{
\arrowlength #3
\put(0,0){\vector(#1,#2){\arrowlength}}%
\advance \arrowlength by-30
\put(0,0){\vector(#1,#2){\arrowlength}}%
}}
 
\def\horsize#1#2{%
\settowidth{\tempdimen}{$#2$}%
#1=\tempdimen
\divide #1 by\unitlength
}
 
\def\vertsize#1#2{%
\settoheight{\tempdimen}{$#2$}%
#1=\tempdimen
\settodepth{\tempdimen}{$#2$}%
\advance #1 by\tempdimen
\divide #1 by\unitlength
}
 
\def\putvector(#1,#2)(#3,#4)#5#6{{%
\ifnum3<\arrowtype
\putdashvector(#1,#2)(#3,#4)#5\arrowtype
\else
\ifnum\arrowtype<-3
\putdashvector(#1,#2)(#3,#4)#5\arrowtype
\else
\xpos=#1
\ypos=#2
\run=#3
\rise=#4
\arrowlength=#5
\ifnum \arrowtype<0
    \ifnum \run=0
        \advance \ypos by-\arrowlength
    \else
        \tempcounta \arrowlength
        \multiply \tempcounta by\rise
        \divide \tempcounta by\run
        \ifnum\run>0
            \advance \xpos by\arrowlength
            \advance \ypos by\tempcounta
        \else
            \advance \xpos by-\arrowlength
            \advance \ypos by-\tempcounta
        \fi
    \fi
    \multiply \arrowtype by-1
    \multiply \rise by-1
    \multiply \run by-1
\fi
\ifcase \arrowtype
\or \put(\xpos,\ypos){\vector(\run,\rise){\arrowlength}}%
\or \put(\xpos,\ypos){\mvector(\run,\rise)\arrowlength}%
\or \put(\xpos,\ypos){\evector(\run,\rise){\arrowlength}}%
\fi\fi\fi
}}
 
\def\putsplitvector(#1,#2)#3#4{
\xpos #1
\ypos #2
\arrowtype #4
\halflength #3
\arrowlength #3
\gap 140
\advance \halflength by-\gap
\divide \halflength by2
\ifnum\arrowtype>0
   \ifcase \arrowtype
   \or \put(\xpos,\ypos){\line(0,-1){\halflength}}%
       \advance\ypos by-\halflength
       \advance\ypos by-\gap
       \put(\xpos,\ypos){\vector(0,-1){\halflength}}%
   \or \put(\xpos,\ypos){\line(0,-1)\halflength}%
       \put(\xpos,\ypos){\vector(0,-1)3}%
       \advance\ypos by-\halflength
       \advance\ypos by-\gap
       \put(\xpos,\ypos){\vector(0,-1){\halflength}}%
   \or \put(\xpos,\ypos){\line(0,-1)\halflength}%
       \advance\ypos by-\halflength
       \advance\ypos by-\gap
       \put(\xpos,\ypos){\evector(0,-1){\halflength}}%
   \fi
\else \arrowtype=-\arrowtype
   \ifcase\arrowtype
   \or \advance \ypos by-\arrowlength
       \put(\xpos,\ypos){\line(0,1){\halflength}}%
       \advance\ypos by\halflength
       \advance\ypos by\gap
       \put(\xpos,\ypos){\vector(0,1){\halflength}}%
   \or \advance \ypos by-\arrowlength
       \put(\xpos,\ypos){\line(0,1)\halflength}%
       \put(\xpos,\ypos){\vector(0,1)3}%
       \advance\ypos by\halflength
       \advance\ypos by\gap
       \put(\xpos,\ypos){\vector(0,1){\halflength}}%
   \or \advance \ypos by-\arrowlength
       \put(\xpos,\ypos){\line(0,1)\halflength}%
       \advance\ypos by\halflength
       \advance\ypos by\gap
       \put(\xpos,\ypos){\evector(0,1){\halflength}}%
   \fi
\fi
}
 
\def\putmorphism(#1)(#2,#3)[#4`#5`#6]#7#8#9{{%
\run #2
\rise #3
\ifnum\rise=0
  \puthmorphism(#1)[#4`#5`#6]{#7}{#8}#9%
\else\ifnum\run=0
  \putvmorphism(#1)[#4`#5`#6]{#7}{#8}#9%
\else
\setpos(#1)%
\arrowlength #7
\arrowtype #8
\ifnum\run=0
\else\ifnum\rise=0
\else
\ifnum\run>0
    \coefa=1
\else
   \coefa=-1
\fi
\ifnum\arrowtype>0
   \coefb=0
   \coefc=-1
\else
   \coefb=\coefa
   \coefc=1
   \arrowtype=-\arrowtype
\fi
\width=2
\multiply \width by\run
\divide \width by\rise
\ifnum \width<0  \width=-\width\fi
\advance\width by60
\if l#9 \width=-\width\fi
\putbox(\xpos,\ypos){#4}
{\multiply \coefa by\arrowlength
\advance\xpos by\coefa
\multiply \coefa by\rise
\divide \coefa by\run
\advance \ypos by\coefa
\putbox(\xpos,\ypos){#5} }%
{\multiply \coefa by\arrowlength
\divide \coefa by2
\advance \xpos by\coefa
\advance \xpos by\width
\multiply \coefa by\rise
\divide \coefa by\run
\advance \ypos by\coefa
\if l#9%
   \putrbox(\xpos,\ypos){#6}%
\else\if r#9%
   \putlbox(\xpos,\ypos){#6}%
\fi\fi }%
{\multiply \rise by-\coefc
\multiply \run by-\coefc
\multiply \coefb by\arrowlength
\advance \xpos by\coefb
\multiply \coefb by\rise
\divide \coefb by\run
\advance \ypos by\coefb
\multiply \coefc by70
\advance \ypos by\coefc
\multiply \coefc by\run
\divide \coefc by\rise
\advance \xpos by\coefc
\multiply \coefa by140
\multiply \coefa by\run
\divide \coefa by\rise
\advance \arrowlength by\coefa
\ifcase\arrowtype
\or \put(\xpos,\ypos){\vector(\run,\rise){\arrowlength}}%
\or \put(\xpos,\ypos){\mvector(\run,\rise){\arrowlength}}%
\or \put(\xpos,\ypos){\evector(\run,\rise){\arrowlength}}%
\fi}\fi\fi\fi\fi}}

\newcount\numbdashes \newcount\lengthdash \newcount\increment
 
\def\howmanydashes{
\numbdashes=\arrowlength \lengthdash=40
\divide\numbdashes by \lengthdash
\lengthdash=\arrowlength
\divide\lengthdash by \numbdashes
\increment=\lengthdash
\multiply\lengthdash by 3
\divide\lengthdash by 5
}
 
\def\putdashvector(#1)(#2,#3)#4#5{%
\ifnum#3=0 \putdashhvector(#1){#4}#5
\else
\ifnum#2=0
\putdashvvector(#1){#4}#5\fi\fi}
 
\def\putdashhvector(#1,#2)#3#4{{%
\arrowlength=#3 \howmanydashes
\multiput(#1,#2)(\increment,0){\numbdashes}%
{\vrule height .4pt width \lengthdash\unitlength}
\arrowtype=#4 \xpos=#1
\ifnum\arrowtype<0 \advance\arrowtype by 7 \fi
\ifcase\arrowtype
\or \advance\xpos by 10
    \put(\xpos,#2){\vector(-1,0){\lengthdash}}
    \advance\xpos by 40
    \put(\xpos,#2){\vector(-1,0){\lengthdash}}
\or \advance \xpos by 10
    \put(\xpos,#2){\vector(-1,0){\lengthdash}}
    \advance\xpos by  \arrowlength
    \advance\xpos by  -50
    \put(\xpos,#2){\vector(-1,0){\lengthdash}}
\or \advance\xpos by 10
    \put(\xpos,#2){\vector(-1,0){\lengthdash}}
\or \advance\xpos by \arrowlength
    \advance\xpos by -\lengthdash
    \put(\xpos,#2){\vector(1,0){\lengthdash}}
\or {\advance\xpos by 10
    \put(\xpos,#2){\vector(1,0){\lengthdash}}}
    \advance\xpos by \arrowlength
    \advance\xpos by -\lengthdash
    \put(\xpos,#2){\vector(1,0){\lengthdash}}
\or \advance\xpos by \arrowlength
    \advance\xpos by -\lengthdash
    \put(\xpos,#2){\vector(1,0){\lengthdash}}
    \advance\xpos by -40
    \put(\xpos,#2){\vector(1,0){\lengthdash}}
   \fi
}}
 
\def\putdashvvector(#1,#2)#3#4{{%
\arrowlength=#3 \howmanydashes
\ypos=#2 \advance\ypos by -\arrowlength
\multiput(#1,#2)(0,\increment){\numbdashes}%
    {\vrule width .4pt height \lengthdash\unitlength}
\arrowtype=#4 \ypos=#2
\ifnum\arrowtype<0 \advance\arrowtype by 7 \fi
\ifcase\arrowtype
\or \advance\ypos by \arrowlength \advance\ypos by -40
    \put(#1,\ypos){\vector(0,1){\lengthdash}}
    \advance\ypos by -40
    \put(#1,\ypos){\vector(0,1){\lengthdash}}
\or \advance\ypos by 10
    \put(#1,\ypos){\vector(0,1){\lengthdash}}
    \advance\ypos by \arrowlength \advance\ypos by -40
    \put(#1,\ypos){\vector(0,1){\lengthdash}}
\or \advance\ypos by \arrowlength \advance\ypos by -40
    \put(#1,\ypos){\vector(0,1){\lengthdash}}
\or \advance\ypos by 10
    \put(#1,\ypos){\vector(0,-1){\lengthdash}}
\or \advance\ypos by 10
    \put(#1,\ypos){\vector(0,-1){\lengthdash}}
    \advance\ypos by \arrowlength \advance\ypos by -40
    \put(#1,\ypos){\vector(0,-1){\lengthdash}}
\or \advance\ypos by 10
    \put(#1,\ypos){\vector(0,-1){\lengthdash}}
    \advance\ypos by 40
    \put(#1,\ypos){\vector(0,-1){\lengthdash}}
\fi
}}
 
\def\puthmorphism(#1,#2)[#3`#4`#5]#6#7#8{{%
\xpos #1
\ypos #2
\width #6
\arrowlength #6
\arrowtype=#7
\putbox(\xpos,\ypos){#3\vphantom{#4}}%
{\advance \xpos by\arrowlength
\putbox(\xpos,\ypos){\vphantom{#3}#4}}%
\horsize{\tempcounta}{#3}%
\horsize{\tempcountb}{#4}%
\divide \tempcounta by2
\divide \tempcountb by2
\advance \tempcounta by30
\advance \tempcountb by30
\advance \xpos by\tempcounta
\advance \arrowlength by-\tempcounta
\advance \arrowlength by-\tempcountb
\putvector(\xpos,\ypos)(1,0)\arrowlength\arrowtype
\divide \arrowlength by2
\advance \xpos by\arrowlength
\vertsize{\tempcounta}{#5}%
\divide\tempcounta by2
\advance \tempcounta by20
\if a#8 %
   \advance \ypos by\tempcounta
   \putbox(\xpos,\ypos){#5}%
\else
   \advance \ypos by-\tempcounta
   \putbox(\xpos,\ypos){#5}%
\fi}}
 
\def\putvmorphism(#1,#2)[#3`#4`#5]#6#7#8{{%
\xpos #1
\ypos #2
\arrowlength #6
\arrowtype #7
\settowidth{\xlen}{$#5$}%
\putbox(\xpos,\ypos){#3}%
{\advance \ypos by-\arrowlength
\putbox(\xpos,\ypos){#4}}%
{\advance\arrowlength by-140
\advance \ypos by-70
\ifdim\xlen>0pt
   \if m#8%
      \putsplitvector(\xpos,\ypos)\arrowlength\arrowtype
   \else
   \putvector(\xpos,\ypos)(0,-1)\arrowlength\arrowtype
   \fi
\else
   \putvector(\xpos,\ypos)(0,-1)\arrowlength\arrowtype
\fi}%
\ifdim\xlen>0pt
   \divide \arrowlength by2
   \advance\ypos by-\arrowlength
   \if l#8%
      \advance \xpos by-40
      \putrbox(\xpos,\ypos){#5}%
   \else\if r#8%
      \advance \xpos by40
      \putlbox(\xpos,\ypos){#5}%
   \else
      \putbox(\xpos,\ypos){#5}%
   \fi\fi
\fi
}}
 
\def\putsquarep<#1>(#2)[#3;#4`#5`#6`#7]{{%
\setsqparms[#1]%
\setpos(#2)%
\settokens`#3`%
\puthmorphism(\xpos,\ypos)[\tokenc`\tokend`{#7}]{\width}{\arrowtyped}b%
\advance\ypos by \height
\puthmorphism(\xpos,\ypos)[\tokena`\tokenb`{#4}]{\width}{\arrowtypea}a%
\putvmorphism(\xpos,\ypos)[``{#5}]{\height}{\arrowtypeb}l%
\advance\xpos by \width
\putvmorphism(\xpos,\ypos)[``{#6}]{\height}{\arrowtypec}r%
}}
 
\def\putsquare{\@ifnextchar <{\putsquarep}{\putsquarep%
   <\arrowtypea`\arrowtypeb`\arrowtypec`\arrowtyped;\width`\height>}}
\def\square{\@ifnextchar< {\squarep}{\squarep
   <\arrowtypea`\arrowtypeb`\arrowtypec`\arrowtyped;\width`\height>}}
\def\squarep<#1>[#2`#3`#4`#5;#6`#7`#8`#9]{{
\setsqparms[#1]
\diagram
\putsquarep<\arrowtypea`\arrowtypeb`\arrowtypec`
\arrowtyped;\width`\height>
(0,0)[#2`#3`#4`{#5};#6`#7`#8`{#9}]
\enddiagram
}}                                                 
\def\putptrianglep<#1>(#2,#3)[#4`#5`#6;#7`#8`#9]{{%
\settriparms[#1]%
\xpos=#2 \ypos=#3
\advance\ypos by \height
\puthmorphism(\xpos,\ypos)[#4`#5`{#7}]{\height}{\arrowtypea}a%
\putvmorphism(\xpos,\ypos)[`#6`{#8}]{\height}{\arrowtypeb}l%
\advance\xpos by\height
\putmorphism(\xpos,\ypos)(-1,-1)[``{#9}]{\height}{\arrowtypec}r%
}}
 
\def\putptriangle{\@ifnextchar <{\putptrianglep}{\putptrianglep
   <\arrowtypea`\arrowtypeb`\arrowtypec;\height>}}
\def\ptriangle{\@ifnextchar <{\ptrianglep}{\ptrianglep
   <\arrowtypea`\arrowtypeb`\arrowtypec;\height>}}
\def\ptrianglep<#1>[#2`#3`#4;#5`#6`#7]{{
\settriparms[#1]
\diagram
\putptrianglep<\arrowtypea`\arrowtypeb`
\arrowtypec;\height>
(0,0)[#2`#3`#4;#5`#6`{#7}]
\enddiagram
}}                                            
 
\def\putqtrianglep<#1>(#2,#3)[#4`#5`#6;#7`#8`#9]{{%
\settriparms[#1]%
\xpos=#2 \ypos=#3
\advance\ypos by\height
\puthmorphism(\xpos,\ypos)[#4`#5`{#7}]{\height}{\arrowtypea}a%
\putmorphism(\xpos,\ypos)(1,-1)[``{#8}]{\height}{\arrowtypeb}l%
\advance\xpos by\height
\putvmorphism(\xpos,\ypos)[`#6`{#9}]{\height}{\arrowtypec}r%
}}
 
\def\putqtriangle{\@ifnextchar <{\putqtrianglep}{\putqtrianglep
   <\arrowtypea`\arrowtypeb`\arrowtypec;\height>}}
\def\qtriangle{\@ifnextchar <{\qtrianglep}{\qtrianglep
   <\arrowtypea`\arrowtypeb`\arrowtypec;\height>}}
\def\qtrianglep<#1>[#2`#3`#4;#5`#6`#7]{{
\settriparms[#1]
\width=\height                                
\diagram
\putqtrianglep<\arrowtypea`\arrowtypeb`
\arrowtypec;\height>
(0,0)[#2`#3`#4;#5`#6`{#7}]
\enddiagram
}}
 
\def\putdtrianglep<#1>(#2,#3)[#4`#5`#6;#7`#8`#9]{{%
\settriparms[#1]%
\xpos=#2 \ypos=#3
\puthmorphism(\xpos,\ypos)[#5`#6`{#9}]{\height}{\arrowtypec}b%
\advance\xpos by \height \advance\ypos by\height
\putmorphism(\xpos,\ypos)(-1,-1)[``{#7}]{\height}{\arrowtypea}l%
\putvmorphism(\xpos,\ypos)[#4``{#8}]{\height}{\arrowtypeb}r%
}}
 
\def\putdtriangle{\@ifnextchar <{\putdtrianglep}{\putdtrianglep
   <\arrowtypea`\arrowtypeb`\arrowtypec;\height>}}
\def\dtriangle{\@ifnextchar <{\dtrianglep}{\dtrianglep
   <\arrowtypea`\arrowtypeb`\arrowtypec;\height>}}
\def\dtrianglep<#1>[#2`#3`#4;#5`#6`#7]{{
\settriparms[#1]
\width=\height                                
\diagram
\putdtrianglep<\arrowtypea`\arrowtypeb`
\arrowtypec;\height>
(0,0)[#2`#3`#4;#5`#6`{#7}]
\enddiagram
}}
 
\def\putbtrianglep<#1>(#2,#3)[#4`#5`#6;#7`#8`#9]{{%
\settriparms[#1]%
\xpos=#2 \ypos=#3
\puthmorphism(\xpos,\ypos)[#5`#6`{#9}]{\height}{\arrowtypec}b%
\advance\ypos by\height
\putmorphism(\xpos,\ypos)(1,-1)[``{#8}]{\height}{\arrowtypeb}r%
\putvmorphism(\xpos,\ypos)[#4``{#7}]{\height}{\arrowtypea}l%
}}
 
\def\putbtriangle{\@ifnextchar <{\putbtrianglep}{\putbtrianglep
   <\arrowtypea`\arrowtypeb`\arrowtypec;\height>}}
\def\btriangle{\@ifnextchar <{\btrianglep}{\btrianglep
   <\arrowtypea`\arrowtypeb`\arrowtypec;\height>}}
\def\btrianglep<#1>[#2`#3`#4;#5`#6`#7]{{
\settriparms[#1]
\width=\height                               
\diagram
\putbtrianglep<\arrowtypea`\arrowtypeb`
\arrowtypec;\height>
(0,0)[#2`#3`#4;#5`#6`{#7}]
\enddiagram
}}
 
\def\putAtrianglep<#1>(#2,#3)[#4`#5`#6;#7`#8`#9]{{%
\settriparms[#1]%
\xpos=#2 \ypos=#3
{\multiply \height by2
\puthmorphism(\xpos,\ypos)[#5`#6`{#9}]{\height}{\arrowtypec}b}%
\advance\xpos by\height \advance\ypos by\height
\putmorphism(\xpos,\ypos)(-1,-1)[#4``{#7}]{\height}{\arrowtypea}l%
\putmorphism(\xpos,\ypos)(1,-1)[``{#8}]{\height}{\arrowtypeb}r%
}}
 
\def\putAtriangle{\@ifnextchar <{\putAtrianglep}{\putAtrianglep
   <\arrowtypea`\arrowtypeb`\arrowtypec;\height>}}
\def\Atriangle{\@ifnextchar <{\Atrianglep}{\Atrianglep
   <\arrowtypea`\arrowtypeb`\arrowtypec;\height>}}
\def\Atrianglep<#1>[#2`#3`#4;#5`#6`#7]{{
\settriparms[#1]
\width=\height                                     
\diagram
\putAtrianglep<\arrowtypea`\arrowtypeb`
\arrowtypec;\height>
(0,0)[#2`#3`#4;#5`#6`{#7}]
\enddiagram
}}
 
\def\putAtrianglepairp<#1>(#2)[#3;#4`#5`#6`#7`#8]{{%
\settripairparms[#1]%
\setpos(#2)%
\settokens`#3`%
\puthmorphism(\xpos,\ypos)[\tokenb`\tokenc`{#7}]{\height}{\arrowtyped}b%
\advance\xpos by\height
\puthmorphism(\xpos,\ypos)[\phantom{\tokenc}`\tokend`{#8}]%
{\height}{\arrowtypee}b%
\advance\ypos by\height
\putmorphism(\xpos,\ypos)(-1,-1)[\tokena``{#4}]{\height}{\arrowtypea}l%
\putvmorphism(\xpos,\ypos)[``{#5}]{\height}{\arrowtypeb}m%
\putmorphism(\xpos,\ypos)(1,-1)[``{#6}]{\height}{\arrowtypec}r%
}}
 
\def\putAtrianglepair{\@ifnextchar <{\putAtrianglepairp}{\putAtrianglepairp%
   <\arrowtypea`\arrowtypeb`\arrowtypec`\arrowtyped`\arrowtypee;\height>}}
\def\Atrianglepair{\@ifnextchar <{\Atrianglepairp}{\Atrianglepairp%
   <\arrowtypea`\arrowtypeb`\arrowtypec`\arrowtyped`\arrowtypee;\height>}}
 
\def\Atrianglepairp<#1>[#2;#3`#4`#5`#6`#7]{{
\settripairparms[#1]
\settokens`#2`
\width=\height                                
\diagram
\putAtrianglepairp                            
<\arrowtypea`\arrowtypeb`\arrowtypec`
\arrowtyped`\arrowtypee;\height>
(0,0)[{#2};#3`#4`#5`#6`{#7}]
\enddiagram
}}
 
\def\putVtrianglep<#1>(#2,#3)[#4`#5`#6;#7`#8`#9]{{%
\settriparms[#1]%
\xpos=#2 \ypos=#3
\advance\ypos by\height
{\multiply\height by2
\puthmorphism(\xpos,\ypos)[#4`#5`{#7}]{\height}{\arrowtypea}a}%
\putmorphism(\xpos,\ypos)(1,-1)[`#6`{#8}]{\height}{\arrowtypeb}l%
\advance\xpos by\height
\advance\xpos by\height
\putmorphism(\xpos,\ypos)(-1,-1)[``{#9}]{\height}{\arrowtypec}r%
}}
 
\def\putVtriangle{\@ifnextchar <{\putVtrianglep}{\putVtrianglep
   <\arrowtypea`\arrowtypeb`\arrowtypec;\height>}}
\def\Vtriangle{\@ifnextchar <{\Vtrianglep}{\Vtrianglep
   <\arrowtypea`\arrowtypeb`\arrowtypec;\height>}}
\def\Vtrianglep<#1>[#2`#3`#4;#5`#6`#7]{{
\settriparms[#1]
\width=\height                                 
\diagram
\putVtrianglep<\arrowtypea`\arrowtypeb`
\arrowtypec;\height>
(0,0)[#2`#3`#4;#5`#6`{#7}]
\enddiagram
}}
 
\def\putVtrianglepairp<#1>(#2)[#3;#4`#5`#6`#7`#8]{{
\settripairparms[#1]%
\setpos(#2)%
\settokens`#3`%
\advance\ypos by\height
\putmorphism(\xpos,\ypos)(1,-1)[`\tokend`{#6}]{\height}{\arrowtypec}l%
\puthmorphism(\xpos,\ypos)[\tokena`\tokenb`{#4}]{\height}{\arrowtypea}a%
\advance\xpos by\height
\puthmorphism(\xpos,\ypos)[\phantom{\tokenb}`\tokenc`{#5}]%
{\height}{\arrowtypeb}a%
\putvmorphism(\xpos,\ypos)[``{#7}]{\height}{\arrowtyped}m%
\advance\xpos by\height
\putmorphism(\xpos,\ypos)(-1,-1)[``{#8}]{\height}{\arrowtypee}r%
}}
 
\def\putVtrianglepair{\@ifnextchar <{\putVtrianglepairp}{\putVtrianglepairp%
    <\arrowtypea`\arrowtypeb`\arrowtypec`\arrowtyped`\arrowtypee;\height>}}
\def\Vtrianglepair{\@ifnextchar <{\Vtrianglepairp}{\Vtrianglepairp%
    <\arrowtypea`\arrowtypeb`\arrowtypec`\arrowtyped`\arrowtypee;\height>}}
\def\Vtrianglepairp<#1>[#2;#3`#4`#5`#6`#7]{{
\settripairparms[#1]
\settokens`#2`
\diagram
\putVtrianglepairp                             
<\arrowtypea`\arrowtypeb`\arrowtypec`
\arrowtyped`\arrowtypee;\height>
(0,0)[{#2};#3`#4`#5`#6`{#7}]
\enddiagram
}}

\def\putCtrianglep<#1>(#2,#3)[#4`#5`#6;#7`#8`#9]{{%
\settriparms[#1]%
\xpos=#2 \ypos=#3
\advance\ypos by\height
\putmorphism(\xpos,\ypos)(1,-1)[``{#9}]{\height}{\arrowtypec}l%
\advance\xpos by\height
\advance\ypos by\height
\putmorphism(\xpos,\ypos)(-1,-1)[#4`#5`{#7}]{\height}{\arrowtypea}l%
{\multiply\height by 2
\putvmorphism(\xpos,\ypos)[`#6`{#8}]{\height}{\arrowtypeb}r}%
}}
 
\def\putCtriangle{\@ifnextchar <{\putCtrianglep}{\putCtrianglep
    <\arrowtypea`\arrowtypeb`\arrowtypec;\height>}}
\def\Ctriangle{\@ifnextchar <{\Ctrianglep}{\Ctrianglep
    <\arrowtypea`\arrowtypeb`\arrowtypec;\height>}}
\def\Ctrianglep<#1>[#2`#3`#4;#5`#6`#7]{{
\settriparms[#1]
\width=\height                               
\diagram
\putCtrianglep<\arrowtypea`\arrowtypeb`
\arrowtypec;\height>
(0,0)[#2`#3`#4;#5`#6`{#7}]
\enddiagram
}}                                           
\def\putDtrianglep<#1>(#2,#3)[#4`#5`#6;#7`#8`#9]{{%
\settriparms[#1]%
\xpos=#2 \ypos=#3
\advance\xpos by\height \advance\ypos by\height
\putmorphism(\xpos,\ypos)(-1,-1)[``{#9}]{\height}{\arrowtypec}r%
\advance\xpos by-\height \advance\ypos by\height
\putmorphism(\xpos,\ypos)(1,-1)[`#5`{#8}]{\height}{\arrowtypeb}r%
{\multiply\height by 2
\putvmorphism(\xpos,\ypos)[#4`#6`{#7}]{\height}{\arrowtypea}l}%
}}
 
\def\putDtriangle{\@ifnextchar <{\putDtrianglep}{\putDtrianglep
    <\arrowtypea`\arrowtypeb`\arrowtypec;\height>}}
\def\Dtriangle{\@ifnextchar <{\Dtrianglep}{\Dtrianglep
   <\arrowtypea`\arrowtypeb`\arrowtypec;\height>}}
\def\Dtrianglep<#1>[#2`#3`#4;#5`#6`#7]{{
\settriparms[#1]
\width=\height                              
\diagram
\putDtrianglep<\arrowtypea`\arrowtypeb`
\arrowtypec;\height>
(0,0)[#2`#3`#4;#5`#6`{#7}]
\enddiagram
}}                                          
\def\setrecparms[#1`#2]{\width=#1 \height=#2}%
 
\def\recursep<#1`#2>[#3;#4`#5`#6`#7`#8]{{\m@th
\width=#1 \height=#2
\settokens`#3`
\settowidth{\tempdimen}{$\tokena$}
\ifdim\tempdimen=0pt
  \savebox{\tempboxa}{\hbox{$\tokenb$}}%
  \savebox{\tempboxb}{\hbox{$\tokend$}}%
  \savebox{\tempboxc}{\hbox{$#6$}}%
\else
  \savebox{\tempboxa}{\hbox{$\hbox{$\tokena$}\times\hbox{$\tokenb$}$}}%
  \savebox{\tempboxb}{\hbox{$\hbox{$\tokena$}\times\hbox{$\tokend$}$}}%
  \savebox{\tempboxc}{\hbox{$\hbox{$\tokena$}\times\hbox{$#6$}$}}%
\fi
\ypos=\height
\divide\ypos by 2
\xpos=\ypos
\advance\xpos by \width
\bfig
\putCtrianglep<-1`1`1;\ypos>(0,0)[`\tokenc`;#5`#6`{#7}]%
\puthmorphism(\ypos,0)[\tokend`\usebox{\tempboxb}`{#8}]{\width}{-1}b%
\puthmorphism(\ypos,\height)[\tokenb`\usebox{\tempboxa}`{#4}]{\width}{-1}a%
\advance\ypos by \width
\putvmorphism(\ypos,\height)[``\usebox{\tempboxc}]{\height}1r%
\efig
}}
 
\def\recurse{\@ifnextchar <{\recursep}{\recursep<\width`\height>}}
 
\def\puttwohmorphisms(#1,#2)[#3`#4;#5`#6]#7#8#9{{%
\puthmorphism(#1,#2)[#3`#4`]{#7}0a
\ypos=#2
\advance\ypos by 20
\puthmorphism(#1,\ypos)[\phantom{#3}`\phantom{#4}`#5]{#7}{#8}a
\advance\ypos by -40
\puthmorphism(#1,\ypos)[\phantom{#3}`\phantom{#4}`#6]{#7}{#9}b
}}
 
\def\puttwovmorphisms(#1,#2)[#3`#4;#5`#6]#7#8#9{{%
\putvmorphism(#1,#2)[#3`#4`]{#7}0a
\xpos=#1
\advance\xpos by -20
\putvmorphism(\xpos,#2)[\phantom{#3}`\phantom{#4}`#5]{#7}{#8}l
\advance\xpos by 40
\putvmorphism(\xpos,#2)[\phantom{#3}`\phantom{#4}`#6]{#7}{#9}r
}}
 
\def\puthcoequalizer(#1)[#2`#3`#4;#5`#6`#7]#8#9{{%
\setpos(#1)%
\puttwohmorphisms(\xpos,\ypos)[#2`#3;#5`#6]{#8}11%
\advance\xpos by #8
\puthmorphism(\xpos,\ypos)[\phantom{#3}`#4`#7]{#8}1{#9}
}}
 
\def\putvcoequalizer(#1)[#2`#3`#4;#5`#6`#7]#8#9{{%
\setpos(#1)%
\puttwovmorphisms(\xpos,\ypos)[#2`#3;#5`#6]{#8}11%
\advance\ypos by -#8
\putvmorphism(\xpos,\ypos)[\phantom{#3}`#4`#7]{#8}1{#9}
}}
 
\def\putthreehmorphisms(#1)[#2`#3;#4`#5`#6]#7(#8)#9{{%
\setpos(#1) \settypes(#8)
\if a#9 %
     \vertsize{\tempcounta}{#5}%
     \vertsize{\tempcountb}{#6}%
     \ifnum \tempcounta<\tempcountb \tempcounta=\tempcountb \fi
\else
     \vertsize{\tempcounta}{#4}%
     \vertsize{\tempcountb}{#5}%
     \ifnum \tempcounta<\tempcountb \tempcounta=\tempcountb \fi
\fi
\advance \tempcounta by 60
\puthmorphism(\xpos,\ypos)[#2`#3`#5]{#7}{\arrowtypeb}{#9}
\advance\ypos by \tempcounta
\puthmorphism(\xpos,\ypos)[\phantom{#2}`\phantom{#3}`#4]{#7}{\arrowtypea}{#9}
\advance\ypos by -\tempcounta \advance\ypos by -\tempcounta
\puthmorphism(\xpos,\ypos)[\phantom{#2}`\phantom{#3}`#6]{#7}{\arrowtypec}{#9}
}}
 
\def\setarrowtoks[#1`#2`#3`#4`#5`#6]{%
\def\toka{#1}
\def\tokb{#2}
\def\tokc{#3}
\def\tokd{#4}
\def\toke{#5}
\def\tokf{#6}
}
\def\hex{\@ifnextchar <{\hexp}{\hexp<1000`400>}}
\def\hexp<#1`#2>[#3`#4`#5`#6`#7`#8;#9]{%
\setarrowtoks[#9]
\yext=#2 \advance \yext by #2
\xext=#1 \advance\xext by \yext
\bfig
\putCtriangle<-1`0`1;#2>(0,0)[`#5`;\tokb``\tokd]
\xext=#1 \yext=#2 \advance \yext by #2
\putsquare<1`0`0`1;\xext`\yext>(#2,0)[#3`#4`#7`#8;\toka```\tokf]
\advance \xext by #2
\putDtriangle<0`1`-1;#2>(\xext,0)[`#6`;`\tokc`\toke]
\efig
}
\makeatother

\newcommand\Span{\operatorname{Span}}
\newcommand\Ker{\operatorname{Ker}}
\newcommand\Image{\operatorname{Im}}

\newcommand\Solk{{\operatorname{Sol}_k}}
\newcommand\Sol{{\operatorname{Sol}}}

\newcommand\C{C^\infty}
\newcommand{\inverse}{^{-1}}
\newcommand\ba\overline

\newcommand\R{{\Bbb R}}
\newcommand\Complex{{\Bbb C}}
\newcommand\RPtwo{{{R}P^2}}

\newcommand{\GGG}{{\cal G}}

\newcommand\CC{{\cal C}}

\newcommand\FF{{\cal F}}
\newcommand\EE{{\cal E}}

\newcommand\II{{{\cal I}}}
\newcommand\JJ{{\cal J}}

\newcommand\PP{{\cal P}}

\newcommand\g{{\frak g}}

\newcommand{\ro}{{r_{\!o}}}

\newcommand{\od}{\stackrel{\mbox {\tiny {def}}}{=}}

\newcommand{\irk}{{{\cal I}_k^{(r)}}}
\newcommand{\mkr}{{M_k^{(r)}}}
\newtheorem{maintheorem}{Theorem}
\newtheorem{thm}{Theorem}[section]
\newtheorem{prop}[thm]{Proposition}
\newtheorem{lemma}[thm]{Lemma}
\newtheorem{cor}[thm]{Corollary}
\newtheorem{dfn}[thm]{Definition}

\newtheorem{rmk}[thm]{Remark}
\newtheorem{ex}{Example}[section]
\title [Orbit Reduction of EDS] {
Orbit Reduction of Exterior Differential Systems, and
group-invariant Variational Problems}
\author{ Vladimir Itskov }
\address{Department of Mathematics\\
University of Minnesota\\
Minneapolis,  MN 55455 USA}
\email{(itskov@@math.umn.edu)}
\date{}
\begin{document}
\newcommand \rcf{{r_{\operatorname{cf}}}}
\newcommand \rpfaff{{r_1}}
\newcommand \rstab{{r_{s}}}
\newcommand \rcomm{{r_o}}
\newcommand \rbasic{{r_b}}
\newcommand \rbold{{{r_\eta}}}
\newcommand \p {{\frak p}}
\newcommand\EL{{\operatorname{EL}}}
\maketitle
\begin{abstract}  For a given PDE system (or an exterior differential system) possessing a Lie group
of internal symmetries the orbit reduction procedure is
introduced. It is proved that the solutions of the  reduced
exterior differential system  are in one-to-one correspondence
with the moduli space of  regular solutions of the original
system.
\par  The isomorphism between  the local characteristic cohomology  of the reduced unconstrained jet space
 and  the Lie algebra cohomology of the symmetry group is established.
\par The group-invariant Euler-Lagrange equations of an invariant variational problem are described as a
composition  of the Euler-Lagrange operators on the reduced jet space and certain other differential operators on
the reduced jet space. The practical algorithm of computing these operators is given.

\end{abstract}
\tableofcontents
\section{Introduction.}
 In this paper we introduce and begin to study the orbit reduction of exterior  differential systems.

Recall, that an exterior differential system  \cite{BCGG} is a
pair $(M,\II)$ where $M$ is a manifold, and $\II\subset\bigwedge
T^*M$ is a graded differentially closed ideal. It is a geometrical
generalization of partial differential equations (in this case
$M$ is a submanifold of the jet space and $\II$ is the contact
ideal). The category of exterior differential systems is bigger then the category of
partial differential equations. It can be shown \footnote{See Example \ref{ex-reduced-eds-R3} in
this paper.} that the category of partial differential equations is not closed under the operation
of the orbit reduction  (to be desribed below).  This gives yet another reason for considering
exterior differential systems .
\par  Let  $\EE=(\Delta,\II)$ be a system of partial differential
equations , or more generally,
 an exterior differential system, invariant under the
 action of a  group  $G$ of internal symmetries .
The action of $G$ on $\EE$ induces a $G$-action  on the space
$\Sol(\EE)$ of the solutions of $\EE$.
Let $\EE^{(r)}=(\Delta^{(r)},\II^{(r)})$ be the $r$-th  order
prolongation of $\EE$. The orbit space   $\Delta^{(r)}/G$
possesses the  structure of an exterior differential system
 induced by the structure of  $\EE^{(r)}$.
\par  It turns out  (see Theorem \ref{thm-moduli})   that for high enough order $r$ of prolongation the
solutions of the reduced system are in one-to-one correspondence
with the moduli space  $\Sol(\EE)/G$  of almost all solutions of
the original system.   This motivates the studying of a
group-invariant PDE system through the  study of its reduced exterior differential
system.
 \par All the results of the present  paper are proved for the case of finite-dimensional Lie group actions.
However we believe  that   the same results  remain valid  for the case of real-analytic actions of
infinite-dimensional groups. The infinite-dimensional group action version of   Theorem \ref{thm-moduli},
also suggests a new approach of studying moduli spaces of any locally  defined geometrical objects.
This will be addressed in some other paper.

\par   The other important  reason for studying the orbit reduction
   is the
  {\it inverse problem of reduction}. By inverse reduction we mean
  the following. Given a certain system of nonlinear PDEs one may
  ask a question whether it is an orbit  reduction of a different system
  of PDEs that has a simpler structure. The questions about the
  solutions of the original system  translate  into
  questions  about the solutions of the "simpler"  system.
  For example it would be
  interesting to identify the class of PDEs which are the orbit
  reduction of an unconstrained jet space.  In this case knowing the inverse reduction
  gives the general solution of the original equations.
 \par  As the very first step towards the  understanding the inverse reduction, we establish the
isomorphism between the local characteristic cohomology of the reduced jet space
and the Lie algebra cohomology of a Lie group of contact
transformations acting on the jet space (see Theorem
\ref{thm-cohomology} in this paper).   This in particular,
implies that in order to  realize a PDE system  having an
infinite-dimensional characteristic cohomology as an orbit
reduction of a jet space one   needs  to consider actions of
infinite-dimensional groups.
\par    The other purpose of the present paper is to understand the group-invariant variational
problems via the orbit reduction. As first observed by Sophus Lie  \cite{Lie},
 the  Euler-Lagrange equations  of  every invariant  variational problem  can be written in terms
of the differential invariants of the group action.    In other words, the Euler-Lagrange equations
of a group-invariant variational problem can be pushed forward to the orbit space.
Surprisingly, up to date there was  no  general understanding of  the meaning of the
pushed forward equations on the orbit space, nor there was  a
 general algorithm of producing the group-invariant Euler-Lagrange equations.

  The reduced jet space has its own calculus of variations  ( for example Euler-Lagrange operators ),
that can be interpreted as a calculus of variations with
constrains imposed by the syzygies of the differential
invariants.  It is well-understood  that all the basic
ingredients of such calculus of variations come from the  edge
complex of the corresponding  Vinogradov spectral sequence
\cite{Cspectral2}.
   We show   (see Theorem \ref{thm-EL} below ) that for every
invariant variational problem the push-forward of the invariant
Euler-Lagrange equations onto the orbit space is a composition of
the Euler-Lagrange operators on the reduced jet space and certain
other differential operators. These other differential operators
come from the  morphism of the two Vinogradov spectral sequences
of the original and the reduced jet spaces. We  also give an
explicit algorithm for  computing these differential operators.
\par Here we would like to note that an  alternative approach based on the Cartan's moving frame
method   is   used by I.~Kogan, and P.~Olver  \cite{KO} for computing invariant Euler-Lagrange
equations.

\section{Reduced Exterior Differential Systems}
\label{sec-reduction}
\subsection{Preliminaries: EDS  and PDEs.}
\label{subsec-EDS}
All the geometrical objects considered in this paper are of class
$C^\infty$ unless stated otherwise. All the considered manifolds
are paracompact.

Let $\II=\cup_{x\in M}\II_x$ be a collection of homogeneous\footnote{
 By saying that the ideal $\II_x$ is homogeneous we mean that in the homogeneous-degree
decomposition $\omega_x=\omega_x^1+..+\omega_x^{\dim(M)}$ of $\omega\in\II_x$ every
homogeneous element $\omega_x^n\in\bigwedge^nT^*M$  belongs to the ideal.}
ideals
$$\II_x=\oplus_{n=1}^{\dim(M)}\II_x^n\subset\bigwedge T_x^*M$$ in
the graded exterior algebra $\bigwedge T^*_xM$. We shall say that
a differential form $\omega\in\Omega^n(M)$ is a section of $\II$
( $\omega\in\Gamma(\II)$ ) if for every $x\in M$,
\,\,$\omega(x)\in \II_x^n$. The sections of $\II$ form a {\it
differential ideal}\, if $\quad d\Gamma(\II)\subset\Gamma(\II)$.
We shall  assume that $\Gamma(\II)$ does not contain any
functions except zero.

\begin{dfn}\label{dfn-EDS} We shall say that $\EE=(M,\II)$
 is an Exterior Differential System
( or  EDS for short) if the  space of sections  of $\II$ is a
differential ideal, and there exists a closed subset $X\subset M$
of zero measure, such that for every connected component $U\subset
(M\setminus X)\,\,\,\,$
$\,\,\,\,\II_{|U}=\cup_{x\in U}\II_x$ is a subbundle of $\bigwedge
T^*U$.
\end{dfn}
  In practice it is convenient to define $\II$ by the generators of $\Gamma(\II)$. We shall say that
$\Gamma(\II)$ is generated by the forms $\omega_1,..,\omega_N$
(the notation is $\II=<\omega_1,..,\omega_n>$) if for every
$\omega\in\Gamma(\II)$ there exist forms
$\alpha_i,\beta_i\in\Omega(M)$ such that \newline
$\omega=\sum_{i=1}^{N}(\omega_i\wedge\alpha_i+d\omega_i\wedge\beta_i)$.
\begin{dfn} A $k$-dimensional solution of $\EE=(M,\II)$ is a connected
$k$-dimensional submanifold
$S\hookrightarrow M$, such that the pullback of $\Gamma(\II)$ to S is zero.
\end{dfn}

\begin{ex}\label{ex-jets1} {\bf Jet spaces.}\textnormal { Let $N$ be a manifold. Consider the $r$-th order  jet space $J_k^rN\to N$
of $k$-dimensional submanifolds together with the standard contact ideal
$\CC^{(r)}\subset \bigwedge T^*J_k^rN$
 (see for example \cite{O2}). For  every $k$-dimensional submanifold $S\stackrel{i_S}\hookrightarrow N$
there is a natural lift $ {j^r{i_S}}:S\hookrightarrow J_k^rN$ such that $\theta\in\Gamma(\CC^{(r)})$ if and only
if the pullback of $\theta$ by $j^ri_S$ is zero for every $k$-dimensional submanifold $S$.
The  lifts  $j^rS={j^r{i_S}}(S)$ are
the solutions of the EDS $(J_k^rN,\CC^{(r)})$. }
\end{ex}

\begin{ex}\label{ex-DE1}{\bf PDE systems.} \textnormal {Let $\Delta\stackrel{\iota}\hookrightarrow  J^r_kN$ be a subbundle
of the jet space $J_k^rN$. This subbundle can be thought of as a system of partial differential equations,
whose solutions are  $k$-dimensional submanifolds  $S\hookrightarrow N$ such that $j^rS\subset\Delta$ .
The lifts $j^rS$ of the solutions of $\Delta$ are  the solutions  of the EDS $\EE=(\Delta,\iota^*\CC^{(r)})$.
}
\end{ex}

Note that since the contact ideals on the jet spaces are always generated by one-forms, not every
EDS is described by the last example. However  the prolongation  \cite{BCGG}
of every EDS is a first-order PDE system.
\par Recall that a ($k$-dimensional)  prolongation \cite{BCGG} of $\EE=(M,\II)$ is an EDS
$\EE_k^{(1)}=(M_k^{(1)},\II_k^{(1)})$, where $M_k^{(1)}$ is a set
of all k-dimensional planes in $TM$ annihilating the ideal $\II$:
$$M_k^{(1)}=\{(P,z)|\;\, z\in M ,\; P\subset T_zM,  \,\, \dim(P)=k,\,\, \mbox{and}\quad {\II_z}_{|P}=0    \}
\stackrel{\iota_1}\hookrightarrow J_k^1M,$$
$$ \II_k^{(1)}=\iota_1^*\CC^{(1)}\,.$$
We shall always assume that $\pi^1:M_k^{(1)}\to M $ is a smooth fiber bundle. Sometimes it
 will mean that we remove some closed subset from $M_k^{(1)}$ to make it smooth.
\par For every $k$-dimensional solution $S\hookrightarrow M$ its lift $j^1S\hookrightarrow J^1_kM$ is a submanifold of
 $M_k^{(1)}$, and is a solution of the prolonged EDS $\EE^{(1)}_k$ .
Conversely, given a solution $S_1\hookrightarrow M^{(1)}_k$ of
the prolonged EDS the natural projection $\pi^1(S_1)\subset M$ is
a solution of the original EDS. However this projection may
"lose" some of its dimension, and may happen not to be a smooth
manifold anymore.

\begin{ex}\label{ex-jets-prolongation} \textnormal {The prolongation of $(J^r_kN,\CC^{(r)})$ is $(J^{r+1}_kN,\CC^{(r+1)})$}\end{ex}
\begin{ex}\label{ex-DE-prolongation}{\bf Prolongation of PDE systems.}
 \textnormal {Consider the Example \ref {ex-DE1}.
Denote by $\pi^r:J_k^rN\to N$ the natural projection. For each
small enough open neighborhood $U\subset J_k^rN$ we may introduce
local coordinates $x^1,..,x^k,u^1,..,u^q$ in
$\pi^r(U)\simeq\R^{k+q}\,$ (this actually means that we
artificially impose a structure of a fiber bundle
$\pi^r(U)\to\R^k$). This choice of the coordinates on the base
$N$ induces  the canonical jet coordinates (see for example
\cite{O2,A1}) $(x^i,u^\alpha,u^\alpha_J)$ (here
$J=(J_1,..,J_{| J|})$ is a multiindex  of length   ${|
J|}\leq r$ ). The contact ideal $\CC^{(r)}$ is generated by the
following 1-forms:
\begin{equation}\label{eq-standard-contact-forms}\{\theta_J^\alpha=du^\alpha_J-u^\alpha_{Ji}dx^i\}_{| J|<r}\end{equation}
Any subbundle   $\Delta\hookrightarrow J_k^rN$ can be represented as a zero level set of functions  $\Delta_\nu\in\C(J_k^rN)$.
 Denote by $\frac{d}{dx^i} : \C(J_k^rN) \to\C(J_k^{r+1}N)$ the total derivatives w.r.t. $x^i$.
The PDE system
$$\label{eq-prolonged-pde}
 \Delta^{(1)}\od \{\;\Delta_\nu=0,\;\, \frac{d}{dx^i}\Delta_\nu=0\;\}\stackrel{\iota '}\hookrightarrow J_k^{r+1}N$$
is called a prolongation of $\Delta$ (see for example \cite{VKL1,O1}). }\end{ex}
\begin{lemma}\label{lemma-pde-prolongation} Let $\Delta\stackrel{\iota}\hookrightarrow J^r_kN$  be a subbundle of $J_k^rN\to N$, then the
prolongation of\quad $\EE=(\Delta,\iota^*\CC^{(r)})\quad$ is $\quad\EE_k^{(1)}=(\Delta^{(1)},{\iota '}^*\CC^{(r+1)})$. \end{lemma}
The proof of this lemma is analogous to the proof for the case $r=1$, given in \cite{BCGG} (Example 6.3, pages 153-154).\bigskip
\par  The prolongation of an EDS can be iterated thus giving a  prolongation tower
$$M\leftarrow M_k^{(1)}\leftarrow M_k^{(2)}\leftarrow \cdots \leftarrow M_k^{(\infty)},$$
where  $M_k^{(\infty)}=\lim_{r\to\infty}\mkr$ is the inverse limit. The last lemma furnishes the following
\begin{cor}\label{cor-eds-prolongation}
 Every prolongation tower of an EDS $\EE=(M,\II)$ can be viewed as a prolongation tower of a first-order
PDE system $M_k^{(1)}\hookrightarrow J_k^1M$. In particular, we have the natural embeddings $\iota_r:\mkr\hookrightarrow J_k^rM$, such
that $\iota_r^*\CC^{(r)}=\irk$.
\end{cor}
\subsection{The reduced EDS.}\label{subsec-orbit-reduction}
Let $\GGG$ be some pseudogroup of local diffeomorphisms acting on
a manifold $M$. We shall say that an EDS \, $\EE=(M,\II)$ is
$\GGG$-invariant if for every $x\in M$, and $g\in\GGG$
\begin{equation}\label{eqn-I-inv} g^*\II_{gx}=\II_x\end{equation}
\begin{rmk} If $\EE=(\Delta,\iota^*\CC^{(r)})$  as in Example \ref{ex-DE1} then the symmetry
group $\GGG$ is usually called a group of internal symmetries \cite{Vin84,AKO}.\end{rmk}

We shall always assume that the orbit space $\ba M\od M/\GGG$ is
again a differentiable manifold (in what follows we shall always
denote the orbit spaces by barred symbols).
 The local coordinates on $\bar M$ may be identified with the $G$-invariant functions on
 $M$.
 The local coordinates on  $\ba{M_k^{(r)}} $ are usually called  the differential
 invariants of order $r$ of the $G$-action.

\begin{prop}\label{prop-reduce}
Let $\EE=(M,\II)$  be a $\GGG$-invariant exterior differential
system, then there exists  an exterior  differential system $\ba \EE=(
\ba M,\ba\II)$, such that $\ba{\II}$ is the maximal ideal satisfying
\begin{equation}\label{eqn-projection} \p^*{\ba{\II}}_{\p(x)}\subset\II_x\quad\quad
\forall x \in M,\end{equation} where $\p:M\to\ba M=M/\GGG$ is
the natural projection .
\end{prop}
\begin{dfn} We shall call $\ba\EE=(\ba M,\ba\II)$ {\it the reduced EDS}.\end{dfn}
We will use in the proof the following simple fact
\begin{lemma}\label{lemma-technical}Let $\EE\to B$ be a vector bundle over a manifold $B$, and
 $\EE_i\subset \EE, \;$ $(i=1,2)$ be two subbundles of $\EE$. Then
there exists a closed subset $X\subset B$ of zero measure such
that for every connected component $U\subset (B\setminus X),
\quad$   $\EE_1\cap\EE_2\to U$ is a subbundle of $\EE\to U$.
\end{lemma}
{\bf Proof of Proposition \ref{prop-reduce}.} For every $x\in M$
define $\ba \II_{\p(x)}\subset\bigwedge T^*_{\p(x)}\ba M $
as the preimage of
\begin{equation}\label{eq-Jx}\,\JJ_x\od\II_x\cap\Image \p_x^*\end{equation}
under $\p_x^*:\bigwedge T^*_{\p(x)}\ba M\to \bigwedge
T_x^*M$.

Let us show that this definition does not depend on the choice of
a particular  $x\in \p\inverse(\p(x))$. Assume
$\p(x_1)=\p(x)$, then there exists a local diffeomorphism
$g\in\GGG$, such that $gx_1=x$, and $g^*\II_x=\II_{x_1}$.
Since $\p(g(x_1))=\p(x_1)$, $$
g^*\p_x^*\ba\II_{\p(x)}=\p_{x_1}^*\ba\II_{\p(x)}\,\,\,,$$
therefore $\quad \p_{x_1}^*\ba\II_{\p(x)}\subset
g^*\II_x=\II_{x_1} $, and $\ba \II$ is
well-defined.\newline It is straightforward to check that
$\Gamma(\ba \II)$ is a differential ideal in $\Omega(\ba M)$. To
show that $\ba \II $ is a subbundle of $\bigwedge T^*\ba M$ over
each connected component of the complement to a closed subset of
zero measure, consider $\JJ=\cup_{x\in M}\JJ_x$ (where $\JJ_x$ is
defined in (\ref{eq-Jx})). By lemma (\ref{lemma-technical}), $\JJ$
is again a subbundle outside a closed subset $X_1\subset M$ of
zero measure . Moreover $X_1$ is $\GGG$-invariant. The latter
implies  that $\ba X\od\p(X_1)$ also has Borel measure zero.
Therefore $\ba \II$ is a subbundle over each connected component
of $\ba M\setminus\ba X$, where $\ba X$ is a closed subset of zero
measure. \qed

\begin{ex}\label{ex-reduced-eds-R3} \textnormal { Consider the action of the abelian group $G=\R^3$ on itself
($M=\R^3$) by translations. Define $\EE=(\R^3,<0>)$. The two-dimensional ($k=2$) $r$-th prolongation of $\EE$ is
the jet space of two-dimensional submanifolds:  $$\EE_2^{(r)}=(J_2^r\R^3,\CC^{(r)}).$$  In order to coordinatize the
orbit spaces $\ba{J_2^{r}\R^3}$, we introduce the coordinates $(x^1,x^2,u)$ in $\R^3$ as well as the standard jet
coordinates $u_J$ in the fibers of $J_2^{r}\R^3$ (here  $J$ is a multiindex ). Note that in fact we restricted our attention to the coordinate chart
$U_r\subset{J_2^{r}\R^3}$ that has a complement of zero Borel measure in ${J_2^{r}\R^3}$.
The orbit space $\ba U_r=U_r/\R^3$ is a Euclidean space with  coordinates $(u_J)_{1\leq| J|\leq r}$.
}\par \textnormal {Denote by \begin{equation}\label{eq-yi}y^i\od u_i,\quad i=1,2\end{equation}   the coordinates on the orbit space
$\ba{J_2^{(1)}\R^3}\simeq\RPtwo$. It is obvious that the
reduced ideal  $\ba{\CC^{(1)}}$ is trivial, thus
$$\ba{\EE_2^{(1)}}=(\RPtwo,<0>).$$
The contact ideal on $J_2^2\R^3$ is generated by the three $\R^3$-invariant 1-forms
\begin{equation}\label{eq-ex-etas1} \eta^1=du-u_1dx^1-u_2dx^2,\end{equation}
\begin{equation}\label{eq-ex-etas2}\eta^2=du_1-u_{11}dx^1-u_{12}dx^2 ,\end{equation}
\begin{equation}\label{eq-ex-etas3} \eta^3=du_2-u_{12}dx^1-u_{22}dx^2. \end{equation}
Let us introduce the coordinates on the fiber of $\ba{J_2^2\R^3}\to\ba{J_2^1\R^3}$:
\begin{equation}\label{eq-va} v^1=u_{11},\;v^2=u_{22},\;v^3=u_{12}.\;\end{equation}
Direct calculations show that the reduced ideal $\ba{\CC^{(2)}}\subset\bigwedge T^*\ba{J_2^2\R^3}$ has no 1-form
 component, however it does have a nontrivial 2-form component, generated by the 2-forms
$\bar\omega_1,\bar\omega_2\in\Omega^2(\ba{J_2^2\R^3})$,
\begin{eqnarray*}\bar\omega_1  =(v^2dy^1-v^3dy^2)\wedge dv^1+(v^1dy^2-v^3dy^1)\wedge dv^3=
\\=(u_{11}u_{22}-u^2_{12})d\eta^2 +(u_{22}\eta^2-u_{12}\eta^3)\wedge du_{11}+(u_{11}\eta^3-u_{12}\eta^2)\wedge du_{12},\end{eqnarray*}
\begin{eqnarray*}\bar\omega_2=(v^2dy^1-v^3dy^2)\wedge dv^3+(v^1dy^2-v^3dy^1)\wedge dv^2=
\\=(u_{11}u_{22}-u^2_{12})d\eta^3 +(u_{22}\eta^2-u_{12}\eta^3)\wedge du_{12}+(u_{11}\eta^3-u_{12}\eta^2)\wedge du_{22}
\end{eqnarray*}
(in fact $\ba{\CC^{(2)}}$ is generated by its 2-form component). Therefore
$$ \ba{\EE_2^{(2)}}\simeq(\RPtwo\times\R^3,<\bar\omega_1,\bar\omega_2>).$$
}\end{ex}

The last example shows that although the original EDS is
generated by 1-forms, the reduced EDS does not necessarily have
the same property. In particular, it may not be a prolongation of
anything. This raises the natural question of whether the
reduction procedure commutes with the prolongation. We address
this question in Theorem \ref{thm-reduction} below.

\begin{dfn}\label{dfn-inf-type}We shall say that an EDS $\EE$ is of infinite type if for every $r>1$
\newline $M_k^{(r)}\to M_k^{(r-1)}$ is a differentiable fiber bundle, and $\dim M_k^{(r)}-\dim M_k^{(r-1)}>0.$
\end{dfn}
\bigskip
 Consider a Lie group
$G$, acting on $M$, and $G$-invariant EDS $\EE=(M,\II)$.
 The action of $G$ on $M$ prolongs to the action on
$M_k^{(r)}$. It is well-known \cite{Ovs,O2} that if $\EE=(M,\II)$
is an infinite-type EDS and the action is effective on open
subsets then the $G$-action is  locally free (i.e. the stabilizers are discrete) almost everywhere on
$\mkr$ for big enough $r$.  The author is not aware of any
example when the action does not eventually become free on high
enough prolongation. Moreover, in the real-analytic category
there are strong indications that every effective action becomes
free on high enough prolongation \cite{free}.
 Throughout this paper we shall adopt the following hypothesis:

\bigskip
 \noindent{\bf The Main Assumptions. }
\begin{tabbing}{\bf 1. }\= $G$ is a Lie group, and the considered EDS $\EE$ is
 of infinite type.\\
{\bf 2.}\> There exists  an integer $\rstab$, and a closed subset
$X\subset
 M_k^{(\rstab)}$ of zero Borel \\ \>measure such that the action of $G$
 is free on $M_k^{(\rstab)}\setminus X$.\\
 {\bf 3.}\> The quotient space $\ba{M_k^{(\rstab)}}=(M_k^{(\rstab)}\setminus X)/G$ is a
 differentiable manifold.
\end{tabbing}

\begin{maintheorem} \label{thm-reduction}
Assume that the main assumptions hold.  Then there exist an
integer $\rcomm\geq 0$ such that for every $r\geq \rcomm$ the
procedure of reduction of $\EE_k^{(r)}$ commutes with the
procedure of prolongation, i.e.
\begin{equation}\label{eq-commute} (\overline {\EE_k^{(r)}})_k^{(1)}=\overline{\EE_k^{(r+1)}}
\end{equation}
\end{maintheorem}
\par It will be shown below (see the proof in the section \ref{sec-syzygies} ) that
$\ro\leq\max(\rstab,\rcf)+2$, where $\rcf $ is the order of
prolongation at which a closed horizontal $G$-invariant coframe
appears (see Lemma \ref{lemma-hor-coframe} below).
\subsection{Moduli space of solutions and the reduced EDS}
Let G be a Lie group acting on M. Let $\EE=(M,\II)$ be a
$G$-invariant  EDS of infinite type. Denote by $\Solk(\EE)$ the
space  of $k$-dimensional solutions of
$\EE$.  We shall say that a solution $S_r\in\Solk(\EE^{(r)}_k)$ is
{\it regular} (the notation is
$S_r\in\operatorname{Sol}_k^{\operatorname{reg}}(\EE^{(r)}_k,G)$ )
if  $S_r$  is transversal to the orbits of the $G$-action on
$M^{(r)}_k$ (clearly then the  lifts  of $S_r$ to  the higher
prolongations are  also regular ).
\par For every $r>0$, and every solution $S_r\in\Solk(\EE^{(r)}_k)$ we may consider
 the  projection $\bar S_r=\p(S_r)\subset \ba{\mkr}$ of $S$ onto the orbit space.
If the solution $S_r$ is regular
then $\bar S_r$ is a $k$-dimensional submanifold, and is a solution of the reduced EDS
$\ba{\EE_k^{(r)}}$. It turns out that on "high  enough" prolongation we can also
lift a solution of $\ba{\EE^{(r)}_k}$ to a regular solution of $\EE^{(r)}_k$  .

\begin{maintheorem}\label{thm-moduli} There exists $\ro>0$  (same as in Theorem \ref{thm-reduction}) such that for every
$r\geq\ro$ the moduli space of  regular solutions of the prolonged
EDS is isomorphic to the  solutions of the reduced EDS:
\begin{equation*}\frac{\operatorname{Sol}_k^{\operatorname{reg}}(\EE^{(r)}_k,G)}{G}\simeq\Solk{\ba{\EE^{(r)}_k}}
\;.\end{equation*} \end{maintheorem}
 The proof is given in section \ref{sec-reconstruct}.
\subsection{Characteristic cohomology of the reduced jet spaces.}

Let a Lie group $G$ act on a manifold M. Assume that the main
assumptions hold with regard to the trivial EDS $\EE=(M,<0>)$. By
virtue of Theorem \ref{thm-reduction} we may regard
$(\ba{J_k^{(\infty)}M},\ba{\CC^{(\infty)}})$ as an infinite
prolongation of $\bar\EE_0=(\ba{J_k^{(\ro)}}M,\ba{\CC^{(\ro)}} )$.
\par  The fact that $\bar\EE_0$ is a
reduction of an unconstrained jet space allows us to know
everything about the solutions of $\bar\EE_0$, since every
solution of $\bar\EE_0$ is an image of a solution of $(J_k^\ro
M,\CC^{(\ro)})$ under the mapping $\p :J_k^\ro M\to
\ba{J_k^{(\ro)}M}$. Therefore it is important to investigate the
conditions  under which a given  EDS $\bar\EE_0$ can  be a
reduction of an unconstrained jet space.

\par It turns out that the local characteristic cohomology of the reduced EDS $\bar\EE_0$ is
isomorphic to the Lie algebra cohomology of the Lie group $G$.

\par Denote by ${\bar\pi^\infty_r}:\ba{J_k^\infty M}\to \ba{J_k^r M}$ the natural projection.
\begin{maintheorem}\label{thm-cohomology} For every open subset
$\hat U\subset\ba{ J_k^\infty M}$, such that  ${\bar\pi^\infty_r}\hat U$ is  contractible for
every $r\geq \ro$,
the characteristic cohomology of $\bar\EE_o$ over $\hat U$ is isomorphic to the Lie algebra
cohomology of $G$ in dimensions less than k:
\begin{equation}\label{eq-iso}
H^t(\Omega_{\operatorname{hor}}(\hat U),\bar d_0)\simeq H^t({\frak g})\quad\quad\forall t<k,
\end{equation}
where  $\Omega_{\operatorname{hor}}^t(\hat U)  \od\Omega^t(\hat
U)/ \Gamma(\ba{\CC^{(\infty)}})$, $\;\bar d_0$ is the horizontal
differential induced on the horizontal forms
$\Omega_{\operatorname{hor}}^t(\hat U)$, \; and $H^t({\frak g})$
is the Lie algebra cohomology of the Lie group $G$.
\end{maintheorem}
The proof as well as the practical algorithm of computing the
basis of  nontrivial conservation laws of  $\bar\EE_0$  is given
in section \ref{sec-proof-cohomology} .

\subsection{Invariant variational problems.}
\newcommand\jkm{{J^\infty_kM}}
Consider an unconstrained infinite jet space $J_k^\infty M$ of $k$-dimensional submanifolds  of a manifold $M$.
 Denote by $(\tilde E_r^{s,t},d_r^{s,t})$ the Vinogradov
spectral sequence \cite{Cspectral2}  corresponding to the
decreasing filtration $$\FF_s\Omega(\jkm)=
\Omega(\jkm)\cap\wedge^{s}\Gamma(\CC^{(\infty)}).$$
 It is well-known \cite{Cspectral2,A1} that  the k-dimensional variational problems on
$M$ can be identified with the space $$\tilde
E_1^{0,k}=\Omega^k(\jkm)/\left(\Gamma(\CC^{(\infty)})+d\Omega^{k-1}(\jkm)\right),$$
and the Euler -Lagrange operator is $d_1^{0,k}:\tilde
E_1^{0,k}\to\tilde E_1^{1,k},$ where the quotient
$$\tilde E^{1,k}_1=\Gamma(\CC^{(\infty)})\cap\Omega^{k+1}(\jkm)/\left(\wedge^2\Gamma(\CC^{(\infty)})+
d\Gamma(\CC^{(\infty)})\right)$$
is a free module over the ring of functions on the infinite jet $J_k^\infty M$.
\par For a given $\lambda\in\Omega^k(J_k^rM)$ one may consider the Euler-Lagrange
system \break  $\EL(\lambda)\hookrightarrow J^{2r}_kM$ defined as
the zero locus of $d_1^{0,k}[\lambda]_1$  (we denote by
$[\lambda]_1$  the equivalence class in $\tilde
E_1^{\cdot,\cdot}$).  If $(x^i,u^\alpha,u_J^\alpha)$ are the
standard jet  coordinates in some open neighborhood of $\jkm$,
$\;dx=dx^1\wedge\cdots\wedge dx^k$, and $\lambda=Ldx+
\Gamma(\CC^{(\infty)})$ is the variational problem then the
Euler-Lagrange system has the form
\begin{equation}\label{eq-EL1} \EL(\lambda)=\{ E_\alpha(L)=0, \quad \alpha=1,..,q=\dim M-k  \},\quad \mbox{where}\end{equation}
$$ E_\alpha(L)=\sum_{|I|=0}^{r} (-1)^{|I|}\frac{d^{|I|}}{dx^I}\frac{\partial}{\partial u^\alpha_I}L,
\quad  $$
\begin{equation}\label{eq-EL2} d\lambda=\sum_{\alpha=1}^{q}E_\alpha(L)
\theta^\alpha\wedge dx+d\Gamma(\CC^{(\infty)})+
\wedge^2\Gamma(\CC^{(\infty)}),\quad \end{equation}
\begin{equation*}\mbox{and   the forms}\;
[\theta^\alpha \wedge dx]_1\; \mbox{ give  the basis in } \;
\tilde E_1^{1,k}.\;\end{equation*}
 (Here $\frac{d}{dx^I}$ are the total
derivatives w.r.t.  multiindex $I$.)
   \par Let a Lie group $G$ act on the manifold M.   Since the  $G$-action on $\Omega(\jkm)$
preserves the contact ideal, it  induces the action on $\tilde
E_1^{s,t}$.
\begin{dfn}We shall say that $\lambda\in\Omega^k(\jkm)$
represents  an {\it invariant variational problem} if
$[\lambda]_1$ is $G$-invariant.\end{dfn}

  It can be shown  ( see Lemma \ref{lemma-hor-coframe} in section \ref{sec-proof1})
 that  for every invariant variational problem $[\lambda]_1$ there exists a differential form
$\bar\lambda=\bar Ldy^1\wedge\cdots\wedge dy^k
\in\Omega^k(\ba{J_k^{\infty}M})$ such that
$[\lambda]_1=[\p^*\bar\lambda]_1$ .   The form $\bar \lambda $ in
its turn defines a variational problem  on $\ba{J_k^{\ro}M}$
(that is a class in $\bar E_1^{0,k}$ of the Vinogradov spectral
sequence of $\ba{J^\infty_kM}$). Therefore it is desirable to
understand $\EL(\lambda)$ in terms of the calculus of
variations on the reduced jet space $\ba{J^\infty_kM}$.
\par It is well-known \cite{Tre,O2} that in every small neighborhood of $\ba{J_k^\infty M}$
there exist functions $(y^1,..,y^k,v^1,..,v^{\bar q})$ such that
any other differential invariant is a function of the $y^i$, and
the  total derivatives
$$v^a_I=\frac{d^{|I|}v^a}{dy^I}$$
 of $v^a$ w.r.t. $y^i$ (see Lemma \ref{lemma-tresse} in section
\ref{sec-syzygies} ).   Moreover as a consequence of Theorem
\ref{thm-reduction} we have
$$(\ba{\jkm},\ba{\CC^{(\infty)}})=(\ba{J_k^\ro M},\ba{\CC^{(\ro)}})^{(\infty)}_k
=(\bar\Delta^{(\infty)},\CC^{(\infty)}),$$  where
$\bar\Delta^{(\infty)}\stackrel{\iota_\infty}\hookrightarrow
J_k^\infty\R^{k+\bar q}$ is an infinite prolongation of a certain
PDE system \break  $\bar\Delta\hookrightarrow  J_k^R\R^{k+\bar
q}\;$
 (see section \ref{sec-syzygies} for more details).

In local coordinates the Euler-Lagrange equations on
$$\ba{\jkm}=\bar\Delta^{(\infty)}\stackrel{\iota_\infty}\hookrightarrow J_k^\infty\R^{k+\bar q}$$
 may be written in the same fashion   as on the unconstrained jet space $J_k^\infty\R^{k+\bar q}.$
 More precisely, given a variational problem
$$\bar\lambda=\bar Ldy^1\wedge\cdots\wedge
dy^k\in\Omega^k(\ba{\jkm}),$$ one can find a function
$L_1(y^i,v^a_I)\in\C({J_k^\infty\R^{k+\bar q}})$ such that its
restriction to $\bar\Delta^{(\infty)}$ is equal to $\bar L\;$
($\iota^*_\infty{L_1}=\bar L$), then
\begin{equation}\label{eq-EL3} d\bar\lambda=\sum_{a=1}^{\bar q}\bar E_a(\bar L)
\bar\theta^a\wedge dy+d\Gamma(\ba{\CC^{(\infty)}})+
\wedge^2\Gamma(\ba{\CC^{(\infty)}}),\quad \end{equation} where
$\bar\theta^a\od dv^a-v^a_idy^i$, $dy=dy^1\wedge\cdots\wedge
dy^k,$ and the expression staying in the place of the
  Euler-Lagrange operator
is defined as
 \begin{equation}\label{eq-reduced-EL}  \bar E_a(\bar L)\od\iota_\infty^*\left(
\sum_{I}(-1)^{|I|}\frac{d^{|I|}}{dy^I}\frac{\partial  L_1}{\partial v^a_I}\right),\end{equation}
and  depends on the particular choice of the function $L_1\;$ \footnote{This happens because
$\bar E_1^{1,k}$ is not necessary  a {\it free} module over the ring of functions on  the reduced
jet space $\ba{\jkm}$. }.
\begin{maintheorem}\label{thm-EL}
 Suppose that the main assumptions hold with regard to the trivial EDS $\;(M,<0>)$.
Then there exist total differential operators on the reduced jet space
$\hat A_\alpha^a:\C(\ba{\jkm})\to \C(\ba{\jkm})$,
$$\hat  A_\alpha^a=\sum_{0\leq|I|\leq\ro-1}A^{aI}_\alpha\frac{d^{|I|}}{dy^I},$$
( here\quad
$A^{aI}_\alpha\in\C(\ba{\jkm}),\quad \alpha=1,..,q=\dim M-k ,\quad  a=1,..,\bar q$ )
such that every invariant variational problem
$[\lambda]_1=[\p^*\bar Ldy^1\wedge\cdots\wedge dy^k]_1$ has its Euler-Lagrange system  as
\begin{equation}\label{eq-EL}  \EL(\lambda)=\p\inverse
\left( \left\{ \;\sum_{a=1}^{\bar q}\hat A_\alpha^a\bar E_a (\bar L)=0,\quad\alpha=1,..,q\;\right\}\right),
\end{equation}
where $\bar E_a $ are the Euler-Lagrange operators (\ref{eq-reduced-EL}) on the reduced jet space.
\end{maintheorem}
\begin{rmk}\label{rmk:nezavisit} \textnormal{ Despite the fact that the Euler-Lagrange
 expressions defined
in the formula (\ref{eq-reduced-EL}) depend on the choice of the Lagrangian $L_1$ the
expression $\;\sum_{a=1}^{\bar q}\hat A_\alpha^a\bar E_a (\bar L)$ does not depend on
this freedom.  } \end{rmk}
 The proof as well as the practical algorithm of computing the operators $\hat A_\alpha^a$
 is given in section \ref{sec-inv-EL}.

\section{Invariant contact forms.} \label{sec-proof1}
 Let $\EE=(M,\II) $ be a $G$-invariant EDS of infinite type.
 Denote by $\EE^{(\infty)}_k=(M^{(\infty)}_k, \II^{(\infty)}_k)$
 its infinite prolongation, and by $\p:\mkr\to\overline{\mkr}$
 denote the natural projection onto the orbit space.

 \begin{lemma}\label{lemma-hor-coframe} There exist \,  $\rcf\geq 0$ , and differential invariants  of the $G$-acton
$y^1,..,y^k\in C^\infty(\ba{M^{(\rcf)}_k})\;$ such
 that $\{\p^*dy^1,..,\p^*dy^k\}$ form a basis \break of    $\;\Omega^1(M^{(\infty)}_k)/\Gamma(\II^{(\infty)}_k)$, and
  $\Omega^t(M_k^{(\infty)})/\Gamma(\II^{(\infty)}_k)\;$ is isomorphic to \break
$\bigwedge^t\Span\{\p^*dy^1,..,\p^*dy^k\}$
 as a $\C(M_k^{(\infty)})$-module.
 \end{lemma}
The proof is completely analogous to the proof of the same fact \cite{O2} about the unconstrained jet spaces
$\;(M,<0>)_k^{(r)}=(J_k^rM,\CC^{(r)})\;$  and therefore omitted.
\par Choosing the differential invariants $y^i$ allows us to introduce
the total differential operators
$\frac{d}{dy^i}:C^\infty(\ba{\mkr})\to C^\infty(\ba{M_k^{(r+1)}})\;$
(here $r\geq \rcf$) in the following way. Denote by
$[\,]_0:\Omega^1(\mkr)\to
\Omega^1(M_k^{(\infty)})/\Gamma(\II_k^{(\infty)})$ the natural
projection to the quotient. Then these operators are defined by
the equality
$$ [\p^*dF]_0=\sum_{i=1}^{k}(\p^*\frac{dF}{dy^i})[\p^*dy^i]_0\,, \quad F\in\C(\ba\mkr) $$
It is easy to show that the functions $\frac{dF}{dy^i}\in
C^\infty(\ba{M_k^{(\infty)}})$ actually belong to
$C^\infty(\ba{M_k^{(r+1)}})$. Note also that the operators $\frac{d}{dy^i}$ commute with each other.

\begin{lemma}\label{lemma-prolonged-ideal-tail} Let $(M_k^{(1)},\II_k^{(1)}) $ be the
 prolongation of an Exterior Differential System $\EE=(M,\II)$. Assume that
$\pi:M_k^{(1)}\to M$ is a surjection. Then at every point
$z\in M_k^{(1)} $ the 1-form component
$\II^{(1)1}=\II_k^{(1)}\cap T^*M_k^{(1)}$ of the prolonged ideal  lies in the
pullback of $T^*M$ by $\pi$:
\begin{equation}\label{eq-tail}\II_z^{(1)1}\subset\pi^*(T^*_{\pi(z)}M),\quad \mbox{and}\end{equation}
\begin{equation}\label{eq-prolonged-dim}\dim(\II_z^{(1)1})=\dim M-k\end{equation}
\end{lemma}
\pf Let us introduce local coordinates $x^1,..,x^k,u^1,..,u^q$ in
some open neighborhood of $M$, and the standard jet coordinates
$u^\alpha_i$ in the fibers of $J^1_kM$. The ideal
$\II^{(1)}_k=\iota^*\CC^{(1)}$  on
$M_k^{(1)}\stackrel{\iota}\hookrightarrow J^{1}_kM$ is
generated by the forms
$$\theta^\alpha=\iota^*(du^\alpha-u^\alpha_idx^i),\quad\quad\alpha=1,..,q=\dim M-k.$$
This proves (\ref{eq-tail}). To see that the forms $\theta^\alpha$ are linearly independent observe
that the forms $\iota^*du^1,..,\iota^*du^q,\iota^*dx^1,..,\iota^*dx^k$ are linearly independent
due to the surjectivity of $\pi$.
\qed
\begin{cor}\label{cor-tail} For  every $r>\rcf$, $\alpha\in\Omega^1(M_k^{(r-1)})$, and
 $\beta\in \Omega^2(M_k^{(r-1)})$
$$ (\pi_{r-1}^r)^*\alpha=\theta +\sum_{i=1}^kf_i\p^*dy^i,$$
\begin{equation} \label{eq-somewhat}(\pi_{r-1}^r)^*\beta=\sum_\alpha\theta^\alpha\wedge \tau_\alpha+
\sum_{1\leq i_1<i_2\leq k} f_{i_1i_2}\p^*(dy^{i_1}\wedge
dy^{i_2}),\end{equation} where
$\theta,\theta^\alpha\in\Gamma(\irk)\cap\Omega^1(\mkr)$,
$\tau_\alpha\in\Omega^1(\mkr)$, and $f_i,f_{i_1i_2}\in\C(\mkr)$.
\end{cor}
\newcommand\pirrr{{\overline\pi_{r-1}^r}}
\par Since the action of $G$ on the fiber bundle
$M_k^{(r-1)}\stackrel{\pi_{r-1}^r}\leftarrow\mkr $ is projectible, there is a surjection $\pirrr$ defined
by the following commutative diagram
$$\square<-1`1`1`-1;800`500>[{M_k^{(r-1)}}`{\mkr}`{\ba{M_k^{(r-1)}}}`\ba{\mkr};\pi^r_{r-1}`\p`\p`\pirrr]
$$
\begin{lemma}\label{lemma-reduced-tail}
Assume that $r\geq\rcf+1$. In an open neighborhood of
$\ba{M_k^{(r-1)}}$ introduce local coordinates
$y^1,..,y^k,v^1,..,v^{\tilde q_r}$ (here the functions $y^i$ are
taken from Lemma \ref{lemma-hor-coframe}). Then the 1-form
component of $\ba{\irk}$ is generated by the forms
\begin{equation}\label{eq-theta-bars}\bar\theta^a=dv^a-v^a_idy^i,\quad\quad a=1,..,\tilde q_r=\dim\ba{M_k^{(r-1)}}-k,\end{equation}
where $v^a_i=\frac{dv^a}{dy^i}\in\C(\ba\mkr)$, and  the dimension of this 1-form component
is equal  to $\tilde q_r=\dim\ba{M_k^{(r-1)}}-k$.
\end{lemma}
\pf Due to the definition of the operators $\frac{d}{dy^i}$, the
forms (\ref{eq-theta-bars}) belong to the reduced ideal. Lemma
\ref{lemma-prolonged-ideal-tail} implies that any form
$\bar\theta\in\Omega^1({\ba\mkr})$ that has nonzero projection
into the quotient $\Omega^1(\ba\mkr )/
(\pirrr)^*\Omega^1(\ba{M_k^{(r-1)}})$ may not belong to the
reduced ideal, thus the forms (\ref{eq-theta-bars})  exhausts the
list of the generators. These forms are linearly independent
 because of the surjectivity of $\pirrr$.
\qed
\begin{lemma}\label{lemma-free-basis} Assume that a Lie Group $G$ acts freely on a manifold
$B$, and has a projectible action on a vector bundle $E\to B$.
Then  every point on the base has an open neighborhood $U\subset
B$ having a basis of $G$-invariant sections in $\Gamma(E_{|U})$.
\end{lemma}
\pf For a given point on the base $B$ choose a small neighborhood
$U\subset B$, such that there exists a right moving frame
\cite{FO2} $\rho:U\to G\,\,$ (i.e.
$\rho(gx)=\rho(x)g^{-1}\,\,\,\,\forall (x,g)\in U\times G$).
Denote by $I(x)=\rho(x)x$ the invariantization map \cite{FO2}. Its
image $L=I(U)$ is a submanifold of U, transversal to the orbits
of the $G$-action on $B$. We may assume that the restriction
vector bundle $E_{|L}\to L $ is trivial. Consider a basis of sections
$\tilde s_1,..,\tilde s_l\in \Gamma(E_{|L}) \,\, (l=\dim(E_x))$
over $L$. Then it is easy to show that the sections
$$s_i(x)=\rho(x)\inverse\tilde s_i(\rho(x)x),\quad\quad i=1,..,\dim(E_x),\quad x\in U,$$
are $G$-invariant and constitute a basis in $\Gamma(E_{|U})$
\qed \newline

Let   $\p:\mkr\to\overline{\mkr}$ denote the natural projection onto the orbit space.
Denote by $\Omega^1(\mkr)^G$ the space of G-invariant differential 1-forms on $\mkr$.

Applying the Lemma \ref{lemma-free-basis} to the first-degree
component of the prolonged ideal $({\irk})^1\subset T^*\mkr$
gives  the following corollary.
\begin{prop}\label{prop-etas}

Suppose that the main assumptions  hold.\newline Denote
$\rbold=\max(\rstab,\rcf)+1$. Then  every point $z\in
M_k^{(\rbold)}\setminus (\pi^\rbold_\rstab)\inverse X$ has  an
open neighborhood $U\subset M_k^{(\rbold)}$, and the $G$-invariant
contact forms
\begin{equation}\label{eq-etas}  \eta^1,\eta^2..,\eta^{\dim G}
\in\Omega^1(U)^G\cap\Gamma(\II_k^{(\rbold)}),
\end{equation}
such that  $$\Span\{\eta^1,\eta^2,..,\eta^{\dim G}\}\simeq\Omega^1(U)^G/ \, \p^*\Omega^1(\p( U)).$$
\end{prop}
\pf Apply Lemma \ref{lemma-free-basis} to the subbundle  $\II^{(\rbold)}_k\cap T^*M_k^{(\rbold)}$.
Then over a certain neighborhood $U\subset M_k^{(\rbold)}$ we have a $G$-invariant basis of
contact forms $\zeta^1,..,\zeta^{q_\rbold}\in\Omega^1(U)^G\cap\Gamma(\II_k^{(\rbold)})\;$,
where $q_\rbold=\dim M_k^{(\rbold-1)}-k \;$  ( see \break Lemma \ref{lemma-prolonged-ideal-tail} ) .
Since $G$ acts freely on $\pi_{\rbold-1}^\rbold(U)$, the difference between the dimensions of the
1-form component of $\II_k^{(\rbold)}$, and the 1-form component of $\ba{\II_k^{(\rbold)}}$ is equal to the dimension of $G$:
$$\dim(\II_k^{(\rbold)})^1-\dim(\ba{\II_k^{(\rbold)}})^1=q_\rbold-\tilde q_\rbold=\dim M_k^{(\rbold-1)}-\dim\ba{M_k^{(\rbold-1)}}=\dim G$$
(here we used Lemma \ref{lemma-reduced-tail} ). Therefore the projection of $\Span\{\zeta^1,..,\zeta^{q_\rbold}\}$  into the quotient
$\Omega^1(U)^G/ \, \p^*\Omega^1(\p( U))$ is surjective. Thus we can find the forms (\ref{eq-etas}) as a
linear combination of the forms $\zeta^j$ with the lifts of some functions on $\ba{M_k^{(\rbold)}}$.
\qed
\begin{cor}\label{cor-etas}For every $r\geq\rbold$, and $z\in M_k^{(\rbold)}\setminus
(\pi^\rbold_\rstab)\inverse X$ each \break
$z'\in({\pi^r_\rbold)}^{-1}(z)$ has an open neighborhood
$U'\subset(\pi^r_\rbold)^{-1}(U)\subset \mkr$, and  contact
forms on the orbit space
$$ \ba\theta^1,..,\ba\theta^{\tilde q_r}\in\Omega^1(\p( U')),$$
such that  the forms $\p^*\bar\theta^a$, and $(\pi^r_\rbold)^*\eta^j$ generate
$\Gamma(\II^{(r)})^G$,  i.e.
 every  $G$-invariant differential form $\omega\in\Gamma(\irk)\cap(\Omega^n(U'))^G$
  satisfies
\begin{equation}\label{eq-freeEDS}
\omega=\sum_{i=1}^{{\dim G}}
((\pi^r_\rbold)^*\eta^i\wedge\alpha_i+(\pi^r_\rbold)^*d\eta^i\wedge\beta_i )+\sum_{a=1}^{\tilde q_r}
(\p^*\ba \theta^a\wedge\gamma_a +  \p^*d\ba\theta^a\wedge\delta_a),
\end{equation} where  the  differential forms \,
 $\alpha_i,\gamma_a\in (\Omega^{n-1}(U'))^G,\, \,\beta_i,\delta_a\in (\Omega^{n-2}(U'))^G\,\:$
are $G$-invariant.
\par Moreover, if $r>\rbold$, then
\begin{equation}\label{eq-twostars}(\pi^r_\rbold)^*d\eta^i=\sum_{j=1}^{\dim
G}(\pi^r_\rbold)^*\eta^j\wedge\alpha_j^i+\sum_{a=1}^{\tilde
q_r}\p^*\bar\theta^a\wedge\zeta^i_a, \end{equation}
where $\quad i=1,..,\dim G ,\quad$  $\alpha_j^i,\zeta^i_a\in \Omega^1(\mkr)^G$.
\end{cor}
\pf  The forms $\bar\theta^a$ are obtained by applying Lemma \ref{lemma-free-basis} to the subbundle
$\II_k^{(r)1}=\irk\cap T^*M_k^{(r)}$, and using the forms $(\pi_\rbold^r)^*\eta^i$ to project the $G$-nvariant
basis of $\Gamma(\II_k^{(r)1})$
to $\p^*(\Omega^1(\ba\mkr))$ thus getting the forms on the orbit space. The formula (\ref{eq-freeEDS}) holds
because $(\Gamma(\irk))^G$ is generated by its 1-form component.
\par If $r>\rbold$, we  can use  formula (\ref{eq-somewhat}),
and write $$(\pi^r_\rbold)^*d\eta^i=\sum_\alpha\theta^\alpha\wedge
\tau_\alpha+\sum_{1\leq i_1<i_2\leq k
}f_{i_1i_2}\p^*(dy^{i_1}\wedge dy^{i_2}) .$$Taking into account that
the right-hand side of the last equation belongs to the ideal
$\Gamma(\irk)$, and using Lemma \ref{lemma-hor-coframe}  we
conclude that $f_{i_1i_2}=0$. Rewriting $\theta^\alpha$ as linear
combination of $(\pi^r_\rbold)^*\eta^i$, and $\bar\theta^a$ gives
the  formula  (\ref{eq-twostars}).
\qed

\begin{ex}\label{ex-R3-etas}\textnormal {Consider the example \ref{ex-reduced-eds-R3}.
Here $\rstab=\rcf=1$, thus $\rbold=2$. The generators
(\ref{eq-standard-contact-forms}) of the contact ideal $\CC^{(r)}$ are already invariant.  However none of them
lives on the orbit space. The forms $\eta^1,\eta^2,\eta^3$ are defined in (\ref{eq-ex-etas1})-(\ref{eq-ex-etas3}).
 Consider $\ba{M_2^{(3)}}=\ba{J^3_2\R^3}$ with the coordinates $(y^i,v^a,v^a_i)$, where
$(y^i,v^a)$ are defined in (\ref{eq-yi}),(\ref{eq-va}), and $v^a_i=\frac{dv^a}{dy^i}$ (here $i=1,2;\;$ $a=1,2,3$).
(Note that there are functional dependencies among $v_i^a$ .)
Then  $\tilde q_3=3$, and the forms $\bar\theta^a$ are given by the formula
$$\bar\theta^a=dv^a-v^a_idy^i. $$
}\end{ex}
\begin{prop}\label{prop-pfaffian}
For every $r\geq\max(\rstab,\rcf)+2,\;\;\;$ $\ba{\irk}$ is
generated by \mbox{1-forms}, i.e. there exist
$\;\ba\theta^1,..,\ba\theta^{\tilde q_r}\in
\Gamma(\ba{\irk})\cap\Omega^1(\ba{\mkr})$ such that  for every
$\;\ba\omega\in\Gamma({\ba\irk}) $
\begin{equation}\label{eq-omegabar}
\ba\omega=\sum_{a=1}^{\tilde q_r}(\ba\theta^a\wedge\ba\alpha_a+d\ba\theta^a\wedge\ba\beta_a)\quad
\mbox{where}\,\, \ba\alpha_a,\ba\beta_a\in \Omega(\ba\mkr).\end{equation}
\end{prop}
\pf It suffices to prove this proposition in a small neighborhood
in $\ba{\mkr}$.
Given the forms (\ref{eq-etas}), define a  $G$-invariant subbundle
\begin{equation}\label{eq-Pspan}\PP=\Span\{\,(\pi_\rbold^r)^*\eta^i\,\}^{\dim G}_{i=1}\subset
T^*\mkr. \end{equation}    There are  $G$-invariant  decompositions
$$ T^*\mkr=\PP\oplus\p^*T^*\ba\mkr,$$
$$\wedge^nT^*\mkr=\bigoplus_{l=0}^n\wedge^l\PP\otimes\wedge^{n-l}\p^*T^*\ba\mkr.$$
The corresponding $G$-invariant projectors \,\,
 $p_n:\bigwedge^nT^*\mkr\to\bigwedge^n\p^*T^*\ba\mkr$
have the following properties:
\begin{equation}\label{equation-p1}
\omega\in\p^*(\Omega^n(\ba\mkr))\quad\Longleftrightarrow\quad p_n\omega= \omega,\,\, \mbox{and}\,\,
\omega \,\,\mbox{is G-invariant};
\end{equation}
\begin{equation}\label{equation-p2}
\omega\;\;\mbox{is G-invariant}\quad\Longrightarrow \quad p_ n\omega=\p^*\ba\omega\;\;
\mbox{for some}\;\;
\ba\omega\in\Omega^n(\ba\mkr);
\end{equation}
\begin{equation}\label{equation-p3}
p_{n_1+n_2}(\omega_1\wedge\omega_2)=(p_{n_1}\omega_1)\wedge(p_{n_2}\omega_2)\quad\quad
\omega_i\in\Omega^{n_i}(\mkr);
\end{equation}
\begin{equation}\label{equation-p4} p_1((\pi^r_\rbold)^*\eta^i)=0\quad\quad \forall i=1,..,\dim G.
 \end{equation}
Now assume that
$\ba\omega\in\Gamma(\ba\irk)\cap\Omega^n(\ba\mkr)$. Using
Corollary \ref{cor-etas}, and the properties
(\ref{equation-p1})-(\ref{equation-p4}) we conclude that\newline
\par
\noindent ${ \p^*\ba\omega=p_n\p^*\ba\omega=}$
\begin{eqnarray*}
 p_n\left(\sum_{i=1}^{{\dim
G}}\left((\pi^r_\rbold)^*\eta^i\wedge\alpha_i+(\pi^r_\rbold)^*d\eta^i\wedge\beta_i
\right)+\sum_{a=1}^{\tilde q_r}
(\p^*\ba \theta^a\wedge\gamma_a +  \p^*d\ba\theta^a\wedge\delta_a)\right)=  \\
 =\sum_{i=1}^{{\dim G}}\left(p_2(\pi_\rbold^r)^*d\eta^i\right)\wedge
p_{n-2}\beta_i+ \sum_{a=1}^{\tilde q_r} (\p^*\ba \theta^a\wedge
p_{n-1}\gamma_a +  \p^*d\ba\theta^a\wedge p_{n-2}\delta_a)=\\
 =\sum_{i=1}^{{\dim G}} p_2\left(\sum_{j=1}^{\dim
G}(\pi^r_\rbold)^*\eta^j\wedge\alpha_j^i+\sum_{a=1}^{\tilde
q_r}\p^*\bar\theta^a\wedge\zeta^i_a\right) \wedge p_{n-2}\beta_i+ \\
+\sum_{a=1}^{\tilde q_r} (\p^*\ba \theta^a\wedge
p_{n-1}\gamma_a +  \p^*d\ba\theta^a\wedge p_{n-2}\delta_a) =\\
=\sum_{a=1}^{\tilde q_r}\left(\p^*\bar\theta^a\wedge
p_{n-1}\left(\sum_{i=1}^{{\dim G}}\zeta^i_a\wedge\beta_i+
\gamma_a\right) + \p^*d\ba\theta^a\wedge p_{n-2}\delta_a\right).
\end{eqnarray*}
This together with the property (\ref{equation-p2}) proves  formula (\ref{eq-omegabar}).
 \qed

\section {Syzygies of differential invariants and the proof of  Theorem \ref{thm-reduction}.}
\label{sec-syzygies}

\par The following lemma originally appeared in the work of
 A.~Tresse \cite{Tre} in the context of what was later  called
the jet spaces (See more recent treatment in \cite{O2}).
 It says that the differential invariants of any order are generated by taking total derivatives of
finitely many differential invariants. The proof for the case of jet bundles is given in
\cite{O2} (Theorem 5.49 page 171). The proof for the case of an EDS of infinite type is
 completely analogous, and therefore omitted.
\begin{lemma}\label{lemma-tresse} Assume that the EDS $\EE$ is of infinite type, then
 for every \break $r\geq\max(\rstab,\rcf)+2$
the differential invariants of order $r$ are obtained by taking the total derivatives
of the invariants of order $r-1$:
\begin{equation*}\forall f\in\C (\ba{M_k^{(r)}}) \qquad
f=f(y^i,v^a, \frac{dv^a}{dy^i})\end{equation*}
where\, $(y^i,v^a)$ are  local coordinates on $\ba{M_k^{(r-1)}}$ ,  $ i=1,..,k \,,\, a=1,..,\tilde q $.
\end{lemma}
We may think of $(y^i,v^a,v^a_i)$ as standard jet coordinates on  $J_k^1\ba{M_k^{(r-1)}}$.
These coordinate functions give the mapping $\iota_1:\ba\mkr\to J_k^1\ba{M_k^{(r-1)}}$.
The last lemma implies that $\iota_1 $ is an injective immersion of an  open conull subset of
$\ba\mkr$. The image $\bar\Delta$ of $\iota_1$ is
a PDE system  that can be described locally as a zero locus of functions
$\ba\Delta_\nu\in\C( J_k^1\ba{M_k^{(r-1)}})$. These functions are sometimes called
{\it syzygies of differential invariants} \cite{O2}.
\par\bigskip
\noindent{\bf Proof of  Theorem \ref{thm-reduction}.} Let
$r\geq\ro=\max(\rstab,\rcf)+2$. In some open neighborhood ${\bar
U}\subset\ba\mkr$ consider the functions $(y^i,v^a,v^a_i)$ as in
Lemma \ref{lemma-tresse}. These coordinates define
 the embedding $\iota_1:{\bar U}\hookrightarrow J^1_k(\bar\pi_{r-1}^r{\bar U})$.
\par Since $\p\inverse({\bar U})\to {\bar U}$  is a principle $G$-bundle we may find another principle $G$-bundle
$\p_1:{\tilde U}\to J^1_k(\bar\pi_{r-1}^r{\bar U})$ together with the embedding
$\tilde\iota_1:\p\inverse({\bar U})\hookrightarrow {\tilde U}$ such that we have the following
commutative diagram:
$$\square[\p\inverse({\bar U})`{\tilde U}`{\bar U}`J^1_k(\bar\pi_{r-1}^r{\bar U});\tilde\iota_1`\p`\p_1`\iota_1]$$
(in fact we find $\tilde U\to J^1_k(\bar\pi_{r-1}^r{\bar U}) $ only over a certain tubular
neighborhood of the image $\bar\Delta$ of $\iota_1$). The image of the embedding $\tilde\iota_1$ is the zero locus of the
pullbacks of the functions $\bar\Delta_\nu$.
\par Consider the $G$ invariant coframe $(\tilde\eta^j,\p_1^*dy^i,\p_1^*dv^a,\p_1^*dv_i^a)$ on $\tilde U$, where the forms
$\tilde\eta^j$ are uniquely determined by the condition
$\tilde\iota_1^*\tilde\eta^j={\pi_\rbold^r}^*\eta^j$. Denote the dual basis of vector fields on
$\tilde U$ by  $(H_j,\frac\partial{\partial y^i},\frac\partial{\partial v^a},\frac\partial{\partial v_i^a})$.

In order to construct the prolongation of $\EE_k^{(r)}$ we may consider the coordinates
$(p^j_i,p^a_i,p^a_{i_1i_2})$ in an open subset of the fiber of $J^1_k\tilde U\to\tilde U$
 such that each $k$-dimensional
$P\subset T\tilde U$ is given as
$$P=\Span_{i=1,..,k}\{ \frac{\partial}{\partial y^i}+p^j_iH_j+p^a_i\frac{\partial}{\partial v^a}+p^a_{il}
\frac{\partial}{\partial v^a_l} \}$$
(here we use the standard summation convention).
\par Due to  Corollary \ref{cor-etas} the ideal $\Gamma(\irk)$
is algebraically generated by the forms  $({\pi_\rbold^r})^*\eta^j$
, $\bar\theta^a$, and $d\bar\theta^a$.  The embedding
$\p\inverse(\bar U)\hookrightarrow\tilde U$ induces the embedding
$J_k^1\p\inverse(\bar U)\hookrightarrow J_k^1\tilde U$, thus we
may view $M_k^{(r+1)}$ as a submanifold of $J^1_k\tilde U$. We
may construct the prolongation of $\EE^{(r)}_k$ as the
prolongation of  the ideal in $\Omega(\tilde U)$
generated by the forms $\tilde\eta^j$,
$\tilde\theta^a=\p^*_1(dv^a-v^a_idy^i)$,
 $d\tilde\theta^a$, and the functions $\p_1^*\bar\Delta_\nu$. Direct calculation shows that
this prolongation
$M_k^{(r+1)}\hookrightarrow J^1_k\tilde U$ is defined by the equations
\begin{equation}\label{eq-def-eq}\p_1^*\bar\Delta_\nu=0,\quad  \p_1^*\frac{d}{dy^i}\bar\Delta_\nu=0,\end{equation}
\begin{equation}\label{eq-def-eq1}p^j_i=0,\quad p^a_i-v^a_i=0,\quad p^a_{i_1i_2}-p^a_{i_2i_1}=0.\end{equation}
The prolonged ideal $\Gamma(\II^{(r+1)}_k)$ is obtained as the restriction of the standard contact
ideal on $J^1_k\tilde U$ to the zero locus of these functions, and is generated by the forms
$\tilde\eta^i, \;\tilde\theta^a=dv^a-v^a_idy^i,\;$ and  $\tilde\theta^a_i=dv^a_i-p^a_{ii_1}dy^{i_1} $.
\par Due to their definition, the functions $p^j_i,p^a_i,p^a_{i_1i_2}\in\C(J^1_k\tilde U)$ are
$G$-invariant, thus the equations (\ref{eq-def-eq}), (\ref{eq-def-eq1}) can be pushed
forward  by $\p_1$, and the reduced EDS $\ba{\EE_k^{(r+1)}}=(\ba{M_k^{(r+1)}},\ba{\II_k^{(r+1)}})$
 is described by the following data:
$$\ba {M_k^{(r+1)}}=\{\bar\Delta_\nu=0,\; \frac{d\bar\Delta_\nu}{dy^i}=0, \; p^a_{i_1i_2}-p^a_{i_2i_1}=0\}
\hookrightarrow \R^{\tilde q _r+k}_{(y^i,v^a,v^a_i,p^a_{i_1i_2})}\;, $$
$$ \ba{\II_k^{(r+1)}}=<\bar\theta^a=dv^a-v^a_idy^i,\;\bar\theta^a_i=dv^a_i-p^a_{ii_1}dy^{i_1} >$$
(here we used Proposition \ref{prop-pfaffian}  to find $\ba{\II_k^{(r+1)}}$).
It is easy to see  (Lemma \ref{lemma-pde-prolongation}) that this is exactly the prolongation of PDE system defined by the syzygies
$\bar\Delta_\nu$, thus
$$\ba{\EE_k^{(r+1)}}=(\bar\Delta,\iota_1^*\CC^{(1)})^{(1)}_k=(\ba\EE_k^{(r)})^{(1)}_k$$.
\qed
\begin{ex}\label{ex-reduced-eds-R3-sysygies}\textnormal {Consider the example
\ref{ex-reduced-eds-R3}. On the space  $\ba{J^2_2\R^3}$ we introduced the local
 coordinates $(y^1,y^2, v^1,v^2,v^3)$ . Since $\rbold=2$, all the higher order differential invariants
are generated by the total derivatives of  $v^a$.  Counting the
dimensions shows that there are two functionally independent
syzygies, namely
\begin{equation}\label{eq-syzygy1}\bar\Delta_1=v^3(v_2^2-v_1^3)+v^1v_1^2-v^2v^3_2=0,\end{equation}
\begin{equation}\label{eq-syzygy2}\bar\Delta_2=v^3(v_1^1-v_2^3)+v^2v_2^1-v^1v^3_1=0.\end{equation}
Theorem \ref{thm-reduction}  implies that for every $r\geq 3$ the
reduced EDS $\ba{\EE^{(r)}}={{(\ba{J_2^r\R^3}},\ba{\CC^{(r)}})}$
is isomorphic\footnote {In order to fit everything into one
coordinate chart  we actually cut off certain  closed subset of
zero measure from $\ba{J_2^{(r)}\R^3}.$}
 to the $(r-3)$-th
prolongation $(\bar\Delta^{(r-3)},\iota_{r-2}^*\CC^{(r-2)})$  of the PDE system
$\bar\Delta=\{\bar\Delta_1=\bar\Delta_2=0\}\stackrel{\iota_1}\hookrightarrow J_2^1\R^5$.}
\end{ex}
\section{Reconstructing the
solutions of the original EDS and the proof of Theorem
 \ref{thm-moduli}}\label{sec-reconstruct}
\newcommand\pii{\p^{-1}(\bar S)}
\begin{dfn}Let $\JJ^1\subset T^*M $ be a subbundle of the cotangent bundle. Denote by
$\JJ\subset\bigwedge T^*M$ the ideal generated by $\JJ^1$. The EDS
$(M,\JJ)$ is called Frobenius if $\quad
d\Gamma(\JJ^1)\subseteq\Gamma(\JJ^1)\bigwedge\Omega^1(M)$.\end{dfn}
In this case the manifold M is foliated by $k=(\dim M-
\dim\JJ^1_x)$-dimensional solutions of $(M,\JJ)$.

Let $\bar S\hookrightarrow \ba\mkr$ be a $k-dimensional$ solution
 of the reduced EDS $\ba{\EE^{(r)}_k}$. Consider
 $\pii\stackrel{i}\hookrightarrow\mkr$. Define $\JJ(\bar S)\od i^*\irk\subset\bigwedge T^*\pii$.
\begin{prop}\label{prop-frobenius} Let $r\geq\max(\rstab,\rcf)+2$, then the exterior differential
system $(\pii,\JJ(\bar S))$ is Frobenius . The solutions of this
EDS are transversal to the orbits of the $G$-action, and form a
foliation of codimension $\dim G$.
\end{prop}
\pf Since $\bar S $ is a solution of the reduced EDS,
$\;i^*\p^*\bar\theta^a=0$. We can apply the mapping $i^*$ to the
both sides of the equation (\ref{eq-twostars}),
 and conclude that
$$d\,i^*(\pi^r_\rbold)^*\eta^i=\sum_{j=1}^{\dim
G}i^*(\pi^r_\rbold)^*\eta^j\wedge i^*\alpha_j^i.$$
 Therefore the ideal $\JJ(\bar S)$ is algebraically generated by the 1-forms
$\{i^*(\pi^r_\rbold)^*\eta^j\}$.
This proves that the EDS $(\pii,\JJ(\bar S))$ is Frobenius.
\par  Let $\bar x \in \ba\mkr$. At every  point $x\in\p\inverse(\bar x)\subset\pii$ the ideal
$\JJ_x(\bar S)$ is generated by its 1-form component
$$\JJ_x^1(\bar S)=\Span\{i^*(\pi^r_\rbold)^*\eta^i(x)\}_{i=1}^{\dim G}.$$
The commutative diagram
$$\square<-1`1`1`-1;800`500>[{\bar S}`{\pii}`{\ba{M_k^{(r)}}}`\mkr;\p_{\bar S}``i`\p]$$
where the horizontal rows are principal $G$-bundles, and the vertical rows are embeddings,
gives the following commutative diagram
$$ \bfig\putsquare<1`0`-1`1;500`400>(-3100,600)[0`{T_{\bar x}^*\bar
S}`0`T^*_{\bar x}{\ba{\mkr}};```]
\putmorphism(-2450,600)(1,0)[``\p^*]{400}1a
\putmorphism(-2450,1000)(1,0)[``\p^*_{\bar S}]{400}1a
\putsquare<1`-1`-1`1;850`400>(-1800,600)[T^*_x\pii`T^*_x\pii\;/\;\p^*T^*_{\bar
x}\bar S` T^*_x\mkr`T^*_x\mkr\;/\;\p^*T^*_{\bar x
}{\ba{\mkr}};`i^*`\simeq`] \put(0,600){\makebox(0,0){0}}
\put(0,1000){\makebox(0,0){0}}
\putmorphism(-450,600)(1,0)[``]{400}1a
\putmorphism(-450,1000)(1,0)[``]{400}1a
\putsquare<0`0`0`0;400`400>(-1500,0)[``\PP_x`;```]
\putsquare<0`0`0`0;400`400>(-1500,1600)[``\JJ_x^1(\bar S)`;```]
\putmorphism(-1500,1500)(0,-1)[``i^*]{1400}{-1}r
\putmorphism(-950,500)(-1,-1)[``\simeq]{500}{-1}r
\putmorphism(-990,1200)(-1,1)[``]{250}{-1}r
\putmorphism(-1890,500)(2,-3)[``]{300}{-1}r
\putmorphism(-1790,1200)(2,3)[``]{150}{-1}r \efig
 $$
 where the horizontal rows are exact, the leftmost
vertical arrow is epimorphic, and $\Ker i^*\subset \Image \p^*$.
It is easy to see that this implies that the rightmost vertical
arrow is a bijection, thus $\dim \JJ_x^1(\bar S)=\dim \PP_x=\dim
G$ (here the subbundle $\PP\subset T^*\mkr$ is defined in
(\ref{eq-Pspan}) ). The transversality of the solutions and the
group orbits follows from the decomposition $T^*_x\pii=\Image
\p^*_{\bar S}\oplus\JJ_x^1(\bar S)$. \qed \newline

\noindent To prove
Theorem \ref{thm-moduli} we need the following simple
\begin{lemma}\label{lemma-foliation} Let $\p:B\to\bar B$ be a principal $G$-bundle.
Assume that there exists a
$G$-invariant Frobenius EDS $(B,\JJ)$ such that $\dim
\JJ^1_x=\dim G$, and the leaves of the foliation defined by
$(B,\JJ)$ are transversal to the fibers of $\p$. Then for every two connected leaves $S_1, S_2$ of
this foliation there exists a group element $g\in G$ such that
$gS_1=S_2$.
\end{lemma}
 {\bf Proof of Theorem \ref{thm-moduli}. }For every
$S_r\in\operatorname{Sol}_k^{\operatorname{reg}}(\EE^{(r)}_k,G)$
the projection $\bar S=\p(S_r)$ depends only on the equivalence
class of $S_r$ in
$\;\displaystyle{\frac{\operatorname{Sol}_k^{\operatorname{reg}}(\EE^{(r)}_k,G)}{G}}\;$,
and is a $k$-dimensional solution of the reduced EDS.
 Given a solution $\bar S\in \Solk(\ba{\mkr})$
consider a $k$-dimensional solution   $S_r\hookrightarrow\p\inverse(\bar
S)\hookrightarrow\mkr$ of the Frobenius EDS $(\pii,\JJ(\bar S))$.
Proposition \ref{prop-frobenius} implies that  $S_r$ is a regular
$k$-dimensional solution of $\EE_k^{(r)}$. Lemma
\ref{lemma-foliation} implies that a different choice of a
solution of $(\pii,\JJ(\bar S))$ lies in the same equivalence
class of the moduli space. This completes the proof.\qed

\begin{rmk}\label{rmk-parallel-transport} \textnormal { For  every solution $\bar S$ of the reduced EDS
the forms $i^*(\pi^r_\rbold)^*\eta^j\in\Omega^1(\p\inverse(\bar S))$  define
a flat connection in the principle bundle $\p\inverse(\bar S)\to\bar S$.  Thus the
reconstruction of a solution $S_r\hookrightarrow M_k^{(r)}$  constitutes
 finding the parallel transport of a point in $\p\inverse(\bar S)$  w.r.t. this flat connection.
In practical terms this means solving a sequence of $k$ systems of ODEs.  } \end{rmk}

\section {Proof of Theorem \ref{thm-cohomology}, and computing the conservation laws of the syzygy equations.}\label{sec-proof-cohomology}
Let G be a Lie group acting on a manifold M. For every open subset $U\subset J^\infty_k M$ consider a
$G$-invariant Vinogradov spectral sequence $(E_r^{s,t},d_r^{s,t})$ corresponding to the differential
ideal  $\Gamma(\CC^{(\infty)})\cap\Omega(U)^G$ in $\Omega(U)^G$.
 Denote by $(\bar E^{s,t}_r,\bar d_r^{s,t})$ the Vinogradov spectral sequence \cite{Cspectral2,BG}  corresponding
 to the differential ideal $\Gamma(\ba{\CC^{(\infty)}})$ in $\Omega(\p(U))$.
The mapping $\p:J^\infty_k M\to\ba{J^\infty_k M}$ induces the morphism of spectral
sequences  $\p^*:\bar E_r^{s,t}\to E^{s,t}_r$.
\begin{lemma}\label{lemma-horizontal-isomorphism} The map $\p^*$ induces an isomorphism of characteristic
 cohomology:
\begin{equation*} \p^*:\bar E_1^{0,t}\stackrel{\simeq}{\to} E_1^{0,t}.\end{equation*}
 \end{lemma}
\pf Since $d\p^*-\p^*d=0$,\; and
$\p^*(\ba{\CC^{(\infty)}})\subset \CC^{(\infty)}$,\; $\p^*$
induces the mapping $\p^*:\bar E_1^{0,t}{\to} E_1^{0,t}$. Due to
Lemma \ref{lemma-hor-coframe} the mapping $\p^*:\bar
E_0^{0,t}{\to} E_0^{0,t}$ is an isomorphism, thus the induced
mapping in cohomology is also an isomorphism.\qed

The following
theorem first was  announced  in  the paper \cite{AP}  by
I.~Anderson, and J.~Pohjanpelto for the special case when
$M=\R^k\times\R^q$, and the action of $G$ is projectable w.r.t.
the fibration $\R^k\times\R^q\to\R^k$. It turns out that both
these assumptions are superfluous.
\begin{thm}\label{thm-exact-rows} For every open subset $U\subset J^\infty_k M$,
for every $s\geq 1$, and  $t\neq k$  the
 corresponding $G$-invariant Vinogradov spectral sequence $E_r^{s,t}$ satisfies
$$ E_1^{s,t}=0$$\end{thm}
The proof is done by constructing a $G$-invariant  variant of Spencer cohomology,
and proving that it  vanishes for the free complex. The complete proof will be given
elsewhere \cite{my-diser}.
\begin{cor}\label{cor-exact-rows}
\begin{equation}\label{eq-cor-cohomology} E_1^{0,t}\simeq H^t(\Omega(U)^G,d),\quad \;0<t<k,\end{equation}
\begin{equation*} E_2^{s,k}\simeq H^{k+s}(\Omega(U)^G,d),\end{equation*}
where $H^t(\Omega(U)^G,d)$ is the $G$-invariant deRham cohomology of $U\subset J^\infty_k M$.
 \end{cor}
{\bf Proof of Theorem \ref{thm-cohomology}} For every contractible $\hat U$ as in the theorem consider
$U=\p\inverse(\hat U)\subset J_k^\infty M$.  Using  Lemma \ref{lemma-horizontal-isomorphism}, and Corollary
\ref{cor-exact-rows}  we  conclude that for every $t<k,\quad$ $\;\bar E_1^{0,t}\simeq H^t(\Omega(U)^G,d)$.
\par  To prove (\ref{eq-iso}), observe that for every $r>\ro$  $\pi^\infty_rU\to\bar\pi^\infty_r\hat U$ is a principal
 $G$-bundle with contractible base, therefore  $H^t(\Omega(\pi^\infty_rU)^G,d)\simeq H^t({\frak g})$.
This implies that $H^t(\Omega(U)^G,d)\simeq H^t({\frak g})$, thus completing the proof.
\qed
\par Now we would like to describe the practical algorithm for computing the representatives
in the characteristic cohomology classes of the syzygy equations.\par\bigskip

\noindent {\bf The practical algorithm.}\newline {\bf 1.} We may identify $\bigwedge\g^*$ with right-invariant differential
 forms on $G$. Therefore the basis
$\tilde\omega_1,..\tilde\omega_N\in\bigoplus_{t=1}^{k-1}H^t(\g)$ gives the closed forms $\hat\omega_l$ in
$\Omega_{\operatorname{right-inv}}(G)$. For a given contractible subset $\bar U\in\ba{J_k^\rstab M}$ choose
a right moving frame \cite{FO2}, i.e. a mapping $\rho:\bar U\to G$, such that $\rho(gz)=\rho(z)g\inverse$.
The pullbacks $\omega_l=(\pi^\infty_\rstab)^*\rho^*\hat\omega_l$ represent the basis in
$\bigoplus_{t=1}^{k-1}H^t(\Omega(U)^G,d)$.\newline
 {\bf 2.} Using the $G$-invariant coframe $\;\{\eta^j,dy^i,\bar\theta^a\}\;$ in
$\;(\pi^\ro_{\ro-1})^*\Omega^1(J_k^{\ro-1}M)$ we may rewrite each
of the forms $\omega_l$ as $\omega_l=\omega_{li_1\cdots
i_{t_l}}dy^{i_1}\wedge\cdots\wedge dy^{i_{t_l}}
+\Gamma(\CC^{(\ro)})$.\break It is easy to see that the function
$\omega_{li_1\cdots i_{t_l}}$ are $G$-invariant, thus we may
consider the forms on the reduced jet space
$\bar\omega_l\od\omega_{li_1\cdots
i_{t_l}}dy^{i_1}\wedge\cdots\wedge
dy^{i_{t_l}}\in\Omega^{t_l}(\ba{J_k^{(\ro)}M})$. These forms
represent the basis of  characteristic cohomology classes of
$\bar\EE_0$, or using different terminology, nontrivial
conservation laws \cite{O1} of the syzygy equations
$\bar\Delta_\nu=0$.
\begin{ex}\label{ex-cons-laws} \textnormal{Let us compute the nontrivial
conservation laws for the syzygy equations
$\bar\Delta=\{\bar\Delta_1=\bar\Delta_2=0\}\stackrel{\iota_1}\hookrightarrow
J_2^1\R^5$ (\ref{eq-syzygy1}-\ref{eq-syzygy2}) in the example
\ref{ex-reduced-eds-R3-sysygies}.  The Lie algebra cohomology of
$\R^3$ is given by the generators $dx^1,dx^2, du$. The moving
frame $\R^3\to\R^3$ is the multiplication by $-1$ , thus the
forms $dx^1,dx^2, du$ represent the basis in
$H^1(\Omega(J_2^\infty\R^3)^G,d)$. Using the forms
$\eta^1,\eta^2,\eta^3\;$ (\ref{eq-ex-etas1})-(\ref{eq-ex-etas3})
we can notice that
$$ dx^1=\frac{1}{v^1v^2-(v^3)^2}(v^2dy^1-v^3dy^2)+\Gamma(\CC^{(2)}),  $$
$$ dx^2=\frac{1}{v^1v^2-(v^3)^2}(v^1dy^2-v^3dy^1)+\Gamma(\CC^{(2)}), $$
$$ du=y^1dx^1+y^2dx^2 +\Gamma(\CC^{(2)}).$$
Therefore the forms
$$\bar\omega_1=\frac{1}{v^1v^2-(v^3)^2}(v^2dy^1-v^3dy^2),$$
$$\bar\omega_2=\frac{1}{v^1v^2-(v^3)^2}(v^1dy^2-v^3dy^1), $$
$$\bar\omega_3=y^1\bar\omega_1+y^2\bar\omega_2    $$
give the basis of nontrivial conservation laws for the syzygy
equations (\ref{eq-syzygy1}-\ref{eq-syzygy2}). }\end{ex}

\section { Invariant Euler-Lagrange equations.}\label{sec-inv-EL}

\par Consider the increasing filtration  $F_1\subset F_2\subset\cdots\subset F_{\ro}=
 \Gamma(\CC^{(\ro)})\cap\Omega^1(J^{\ro}_kM)$, where
$$F_r\od(\pi^\ro_r)^*\Gamma(\CC^{(r)})\cap\Omega^1(J_k^\ro M).$$
Outside of a certain set of zero measure $F_r$ is a space of
sections of a certain subbundle of $T^*J^\ro_kM$. For each of
these subbundles we can apply Lemma \ref{lemma-free-basis} (in a
small neighborhood of every  point), and find  $G$-invariant
contact forms $\tilde\eta^\alpha_I\in
 \Omega^1(J_k^\ro M)^G\cap\Gamma(\CC^{(\ro)})$ such that for each $r\leq\ro\;$
 $\Span\{\tilde\eta^\alpha_I\}_{|I|\leq r} ^{\alpha=1,..,q}$ is a basis of $F_{r+1}$.
\par Denote by $dy\od dy^1\wedge\cdots\wedge dy^k$ the horizontal  volume on the reduced
jet space.  The following lemma gives the group-invariant version of the integration by parts
used in the  deducing  the Euler-Lagrange equations.
\begin{lemma}\label{lemma-parts} For every $\tilde\eta^\alpha_I$, $|I|>0$ there
exist invariant total differential operators
$\hat T^{\alpha I'}_{I\alpha'}=T^{\alpha I'i}_{I\alpha'}\frac{d}{dy^i}+
T^{\alpha I'}_{I\alpha'}$
 (here $T^{\alpha I'i}_{I\alpha'},T^{\alpha I'}_{I\alpha'}\in\C(\ba{J_k^\ro M})$ ),
such  that for every $\bar f\in\C(\ba{J^\infty_k M})$
\begin{equation} \label{eq:T}(\p^*\bar f)[\tilde\eta^\alpha_I\wedge\p^* dy]_0=\sum_{|I'|<|I|;\, \alpha'=1,..,q}
(\p^*\hat T^{\alpha I'}_{I\alpha'}\bar f )[\tilde\eta^{\alpha'}_{I'}\wedge
\p^*dy]_0+ d_0^{1,k-1}[\chi]_0 \end{equation}
for some  $\chi \in (F_{|I|})^G\wedge\p^*(\C(\ba\jkm)dy) $.
\end{lemma}
The proof is based on the same  fact about the (noninvariant ) standard contact forms $\theta^\alpha_I$
 (\ref{eq-standard-contact-forms}).
\begin{cor}\label{cor-parts}  Let $\bar\theta^a\in\Gamma(\ba{\CC^{\ro}})$ be
the generating  1-forms of the reduced ideal $\ba{\CC^{(\ro)}}$. Then there exist total
differential operators
$$\hat A_\alpha^a:\C(\ba{J_k^rM})\to\C(\ba{J_k^{r+\ro-1}M})\,$$
(here  $r\geq\ro$ ),
$$  \hat A_\alpha^a =\sum_{0\leq|I|\leq \ro-1}A_{\alpha I}^a\frac{d^{|I|}}{dy^I}, \qquad
A_{\alpha I}^a\in\C(J^\infty_kM)$$
such that for every $\bar f\in \C(J_k^\infty M)$
\begin{equation} \label{eq:A}\p^* [\bar\theta^a\wedge\bar fdy]_0=\sum_{\alpha=1}^q
(\p^*\hat A_\alpha^a\bar f)[\tilde\eta^\alpha\wedge\p^*dy]_0+d_0^{1,k-1}[\chi]_0,\end{equation}
for some $[\chi]_0\in E_0^{1,k-1}$, where $\{\tilde\eta^\alpha\}_{\alpha=1,..,q}$  are the basis of forms in
$F_1$.
\end{cor}
\noindent {\bf Proof of  Theorem \ref{thm-EL}.} Let $[\lambda]_1$ be a $G$-invariant
variational problem, then there exists $\bar\lambda=\bar Ldy\in\Omega^k(\ba\jkm)$ such
that $[\p^*\bar\lambda]_1=[\lambda]_1$. Using the above corollary we conclude that
\begin{equation}\label{eq:el3} d_1^{0,k}[\lambda]_1 =d_1^{0,k}[\p^*\bar L dy]_1=\p^*d_1^{0,k}[\bar Ld\lambda]_1
=\p^*[\sum_{a=1}^{\bar q}\bar E_a(\bar L)\bar\theta^a\wedge dy]_1=\end{equation}
\begin{equation*}=\sum_{\alpha=1}^q
\left(\p^*\sum_{a=1}^{\bar q}\hat A_\alpha^a(\bar E_a(\bar L) )\right)[\tilde\eta^\alpha\wedge\p^*dy]_1
=\left(\p^*\hat A_\alpha^a\bar E_a(\bar L)\right)\frac{dy}{dx}c_{\alpha'}^\alpha[\theta^{\alpha'}  \wedge dx ]_1
\end{equation*}
where $\theta^\alpha \in\Gamma(\CC^{(1)})$ are the standard contact forms corresponding to the
choice of local  coordinates $(x^i,u^\alpha)$ on $M$,
 $\;\tilde\eta^\alpha=c^\alpha_{\alpha'}\theta^{\alpha'}$,
and $\p^*dy=\frac{dy}{dx}dx$. Since the matrix  $(\frac{dy}{dx}c_{\alpha'}^\alpha)$ is nondegenerate
the formulas (\ref{eq:el3}), and (\ref{eq-EL1}-\ref{eq-EL2}) imply  (\ref{eq-EL}).
\qed\newline

\noindent {\bf Remark.} Since the functions $(\frac{dy}{dx}c_{\alpha'}^\alpha)$ depend only on the
choice of the horizontal volumes, and the basis of contact forms , the  equality (\ref{eq:el3})
implies that for every $\alpha=1,..,q$ the function
$\sum_{a=1}^{\bar q}\hat A_\alpha^a(\bar E_a(\bar L) )$
does  not depend on the choice of the Lagrangian $L_1$ used in the definition
(  formula (\ref{eq-reduced-EL} )  ) of  $\bar E_a(\bar L)$.\newline

\noindent Now we would like to describe the practical algorithm of computing the \newline
operators $\hat  A^a_\alpha$.
\bigskip

\noindent  {\bf The practical algorithm.}\newline {\bf 1.}  We
can compute the forms $\tilde\eta_I^\alpha$ by applying the
moving frame construction (described in the proof of Lemma
\ref{lemma-free-basis} )  consecutively to each of the subbundles
$(\pi_r^\ro)^*\CC^{(r)1}\subset T^*J^\ro_kM$.\newline {\bf 2.}
For every $r$, $\;0\leq r<\ro$  consider the system of equations
$$(I'\!,\alpha'\!,i')\quad
d_0^{1,k-1}[\tilde\eta^{\alpha'}_{I'}\wedge \nu_i]_0=
\!\!\!\sum_{|I|=r+1}\p^*(\bar f ^{\alpha' I}_{I'\alpha})[\tilde\eta^{\alpha}_{I}\wedge dy]_0
+\sum_{|I''|\leq r}\p^*(\bar f ^{\alpha' I''}_{I'\alpha''})[\tilde\eta^{\alpha''}_{I''}\wedge dy]_0$$
indexed by  the triples  $(I',\alpha',i')\quad $ such that $|I'|=r$, $i=1,..,k$, $\alpha=1,..,q\;$\newline
( here $\nu_i=\p^*(dy^1\wedge\cdots\wedge dy^{i-1}\wedge dy^{i+1}\wedge\cdots\wedge dy^k)$,
and  $\bar f^{\cdot\cdot}_{\cdot\cdot}\in\C(\ba\jkm)$ ).\newline
Due to Lemma \ref{lemma-parts}
we  can always find a solution  $\{\;[\tilde\eta^{\alpha}_{I}\wedge dy]_0\;\}_{|I|=r+1}^{\alpha=1,..,q}$ of this linear overdetermined
 (if $k>1$ ) system of equations and  then using the Leibniz rule compute the operators
$\hat T_{I\alpha'}^{\alpha I'}$  (\ref{eq:T}).\newline
{\bf 3.} We can rewrite the forms $\p^*\bar\theta^a\in (F_\ro)^G$ as a linear combination (over the ring
$\C(\ba\jkm)$ ) of the forms $\tilde\eta^\alpha_I$. Consecutively using the formula (\ref{eq:T}) we
 obtain the operators $\hat A_\alpha^a$  (\ref{eq:A}) .
\begin{ex}\label{ex-SE2} \textnormal {  Consider the ( nonprojectable ) action of the group of
Euclidean motions  $G=\operatorname{SE(2)}$ on $M=\R^2$. Introduce
the standard jet coordinates $(x,u,u_1,u_2,..)$ on
$J_1^\infty\R^2$. The Euclidean curvature
$\kappa=u_2(1+u_1^2)^{-3/2}$, and its derivative with respect to
the arclength $\kappa_s=u_3(1+u_1^2)^{-2}-3u_1u_2^2(1+u_1^2)^{-3}$
provide the local coordinates $y=\kappa, v=\kappa_s $ on the
reduced jet space $\ba{J_1^3\R^2}$. Here $\ro=4$, and  the reduced
EDS $ \ba{\EE^{(4)}}=(\ba{J_1^4\R^2},\ba{\CC^{(4)}})$ is
isomorphic to the first jet space of curves: $
\ba{\EE^{(4)}}=(J^1_1\R^2_{(y,v)},<dv-v_1dy>)$.  Thus the reduced
infinite jet space is again the infinite jet space of curves
\footnote {In fact it is true for any group action on $\R^2$.} in
$\R^2$.  In particular the reduced Euler-Lagrange operators  (\ref{eq-reduced-EL})
coinside with the usual  ones in
$J^\infty_1\R^2.$ }  \par \textnormal {  Let $(c_1,c_2,\phi)$ be the
coordinates on the group $\operatorname{SE(2)}$  such that the action on $\R^2=\Complex$ is given by the formula
$$(c_1,c_2,\phi)(x+iu)=e^{i\phi}(x+iu)+c_1+ic_2.$$
We can use the right moving frame
$\rho:J_1^1\R^2\to\operatorname{SE(2)}$,
$$\rho(x,u,u_1)=(c_1,c_2,\phi)=(-\frac{x+uu_1}{\sqrt{1+u_1^2}},\frac{xu_1-u}{\sqrt{1+u_1^2}},-\tan\inverse u_1)$$
to pull back the right Maurer-Cartan forms
$$\mu^1=d\phi,\;\;\mu^2=dc_1+c_2d\phi,\;\;\mu^3=dc_2-c_1d\phi,$$ and
obtain the basis $$\zeta^i\od\rho^*\mu^i\quad\quad i=1,2,3$$ of
invariant 1-forms in
$\Omega^1(J_1^\infty\R^2)^G/\p^*\Omega^1(\ba{J_1^\infty\R^2})$.
Using the procedure given in the proof of Lemma \ref{lemma-free-basis}
we obtain the filtered basis in $\Gamma(\CC^{(4)})\cap\Omega^1(J_1^4\R^2)^G$:
$$\tilde\eta_0=\zeta^3, \quad \tilde\eta_1=y\zeta^2-\zeta^1,\quad \tilde\eta_2=dy+v\zeta^2,
\quad \tilde\eta_3=\p^*\bar\theta_0$$
(here $y=\kappa$,  $\;v=\kappa_s$, and $\;\bar\theta_0=dv-v_1dy$ ). The table of horizontal
differentiation
$$d_0^{1,0}[\tilde\eta_0]_0=\left[\frac 1v\,\tilde\eta_1\wedge dy\right]_0, $$
$$d_0^{1,0}[\tilde\eta_1]_0=\left[\frac{y^2}{v}\,dy\wedge\tilde\eta_0+\frac 1v dy\wedge\tilde\eta_2 \right]_0 ,$$
$$d_0^{1,0}[\tilde\eta_2]_0=\left[\frac{v_1}v\,dy\wedge\tilde\eta_2+ydy\wedge\tilde\eta_0+\frac{1}{v}dy\wedge\p^*\bar\theta_0 \right]_0 $$
allows us to compute the operators $\hat T$ in the formula (\ref{eq:T}):
$$ (\p^*\bar f)[\p^*(\bar \theta_0\wedge dy)]_1=-\p^*(2v_1\bar f+v\frac {d\bar f}{dy})[\tilde\eta_2\wedge\p^*dy]_1-
\p^*(vy\bar f)[\tilde\eta_0\wedge\p^*dy]_1\,,$$
$$ (\p^*\bar f)[\tilde\eta_2\wedge\p^* dy]_1=-\p^*(v_1\bar f+v\frac {d\bar f}{dy})[\tilde\eta_1\wedge\p^*dy]_1-
\p^*(y^2\bar f)[\tilde\eta_0\wedge\p^*dy]_1\,,$$
$$ (\p^*\bar f)[\tilde\eta_1\wedge\p^* dy]_1=\p^*(v_1\bar f+v\frac{ d\bar f}{dy})[\tilde\eta_0\wedge\p^*dy]_1\,,$$
therefore the operator $\hat A $
 is computed by composing the operators $\hat T$:
$$\hat A=\left((v_1+v\frac d{dy})^2+y^2\right)(2v_1+v\frac d{dy})-vy,$$
and  the Euler-Lagrange system of every invariant variational problem\newline
 $\lambda=Ldx=\bar L dy +\Gamma(\CC^{(\infty)})\;$ is the lift of the  equation  $\hat A
(\sum_{I=0}^\infty(-\frac d{dy})^I\frac{\partial \bar L}{\partial v_I})=0$
on the reduced jet space.
}
\end{ex}

\bigskip
\noindent {\bf Acknowledgments:} The author would like to thank Ian Anderson, Mark Fels,
Irina Kogan, and Peter Olver for interesting   discussions on  the subject.
 \bibliographystyle{plain}
\bibliography{mybib}
\end{document}